\documentclass[12 pt]{article}
\usepackage{amssymb,amsmath,amsfonts,amsthm}
\DeclareMathOperator{\sign}{sign}
\DeclareMathOperator{\BV}{BV}
\numberwithin{equation}{section}
\newcommand{\tq}{\,:\,}

\newcommand{\Monm}{\mathcal{F}}
\newcommand{\tMonm}{\widetilde{\mathcal{F}}}
\DeclareMathOperator{\Mon}{Mon}
\DeclareMathOperator{\tMon}{\widetilde{Mon}}

\newcommand\BBP{{\mathbb {P}}}
\newcommand\p{{\mathbb {P}}}

\newcommand\BBN{{\mathbb {N}}}

\newcommand\bkE{{\mathbb {E}}}
\newcommand\F{{\mathcal {F}}}
\newcommand\N{{\mathbb {N}}}
\newtheorem {Lemma}{Lemma}[section]
\newtheorem {Theorem}{Theorem}[section]
\newtheorem {Proposition}{Proposition}[section]
\newtheorem {Corollary}{Corollary}[section]
\newtheorem{definition}{Definition}[section]

\theoremstyle{definition}
\newtheorem{Definition}{Definition}[section]
\newtheorem{Notation}{Notation}[section]
\newtheorem{Remark}{Remark}[section]

\newcommand\I{{ 1\hspace{-1,2mm}{\mathrm I}}}

\newcommand\ssk{\smallskip}

\newcommand\eps{\varepsilon}
\newcommand\beq{\begin{equation}}
\newcommand\eeq{\end{equation}}

\def\ssk{\smallskip}

\def\cov{\mathop{\rm Cov}\limits}

\def\X{(X_i)_{i\in{\mathbb N}}}
\def\XZ{(X_i)_{i\in{\mathbb Z}}}
\def\N{(N_i)_{i\in{\mathbb Z}}}

\def\eps{\varepsilon}\hfuzz=5pt

\begin{document}
\begin{center} {\bf \Large Strong  approximation of partial sums under dependence conditions with application to dynamical systems}
\vskip15pt

March 15, 2011

Florence Merlev\`{e}de $^{a}$ {\it
and\/} Emmanuel Rio   $^{b}$
\end{center}
\vskip10pt
$^a$ Universit\'e Paris Est-Marne la
Vall\'ee, LAMA and C.N.R.S UMR 8050, 5 Boulevard Descartes, 77454 Marne La
Vall\'ee, FRANCE. E-mail: florence.merlevede@univ-mlv.fr \\ \\
$^b$ Universit\'e de Versailles, Laboratoire de math\'ematiques, UMR
8100 CNRS, B\^atiment Fermat, 45 Avenue des Etats-Unis, 78035
Versailles, FRANCE. E-mail: rio@math.uvsq.fr \vskip10pt

{\it Key words}: almost sure invariance principle, strong approximations, weak dependence, strong mixing, intermittent maps, dynamical systems, Markov chains.\vskip5pt

{\it Mathematical Subject Classification} (2000): 60F17, 37E05.
\begin{center}
{\bf Abstract}\vskip10pt
\end{center}

In this paper, we obtain precise rates of convergence in the strong invariance principle 
for stationary sequences of  real-valued random variables satisfying weak dependence conditions including strong mixing in the sense of Rosenblatt (1956) 
as a special case.  Applications to unbounded functions of intermittent maps are given.

\section{Introduction}
The almost sure invariance principle is a powerful tool in both  probability and statistics. It says that the partial sums of  random variables can be approximated by those of independent Gaussian random variables, and that the approximation error between the trajectories of the two processes is negligible compared to their size. More precisely, when $(X_i)_{i \geq 1}$ is a sequence of  i.i.d. centered  real valued random variables with a finite second moment, a sequence $(Z_i)_{i \geq 1}$ of i.i.d. centered Gaussian variables may be constructed is such a way that 
\beq \label{strassen}
\sup_{1 \leq k \leq n} | \sum_{i=1}^k (X_i - Z_i ) | = o(a_n) \text{ almost surely},
\eeq
where $(a_n)_{n \geq 1}$ is a nondecreasing sequence of positive reals tending to infinity. The first result of this type is due to Strassen (1964) who obtained (\ref{strassen}) with $a_n = (n\log \log n)^{1/2}$. To get smaller $(a_n)$ additional information on the moments of $X_1$ is necessary. If ${\mathbb E}|X_1|^p < \infty$ for $p $ in $]2, 4[$, by using the Skorohod embedding theorem, Breiman (1967) showed that (\ref{strassen}) holds with $a_n = n^{1/p} (\log n)^{1/2}$. He also proved that $a_n =n^{1/p}$ cannot be improved under the $p$-th moment assumption for any $ p >2$. The Breiman paper highlights the fact that there is a gap between the direct result and its converse when using the Skorohod embedding. This gap was later filled by Koml\'os, Major and Tusn\'ady (1976) for $p>3$ and by Major (1976) for $p$ in $]2,3]$: they obtained (\ref{strassen}) with $a_n = n^{1/p}$ as soon as ${\mathbb E}|X_1|^p < \infty$  for any $p>2$, using an explicit construction of the Gaussian random variables, based
on quantile transformations.

There has been a great deal of work to extend these results to dependent sequences: see for instance Philipp and Stout (1975), Berkes and Philipp (1979), Dabrowski (1982), Bradley (1983), Shao (1993), Eberlein (1986),  Wu (2007), Zhao and Woodroofe (2008) among others,  for extensions of (\ref{strassen})  under various dependence conditions. 

In this paper, we are interested in the case of strictly stationary strongly mixing sequences.  Recall that the strong mixing coefficient of Rosenblatt (1956) between two $\sigma$-algebras ${\mathcal F}$
and ${\mathcal G}$ is defined by
$$  \alpha({\mathcal F}, {\mathcal G})= \sup_{A \in {\mathcal F}, B \in {\mathcal G}}|{\mathbb P}(A \cap B)-{\mathbb P}(A){\mathbb P}(B)| \, .
$$
For a strictly stationary sequence $(X_i)_{i \in {\mathbb Z}}$ of real valued random variables, and the $\sigma$-algebra ${\mathcal F}_0=\sigma (X_i, i \leq 0)$ and ${\mathcal G}_n
=  \sigma (X_i, i \geq n)$, define then
\begin{equation}\label{defalpharosen}
 \alpha(0) = 1 \text{ and } \alpha(n)= 2\alpha({\mathcal F}_0,{\mathcal G}_n) \text{ for $n>0$} \, .
\end{equation}
Concerning the extension of (\ref{strassen}) in the strong mixing setting, Rio (1995-a) proved the following: assume that 
\beq \label{quantrio}
\sum_{k=0}^\infty \int_0^{\alpha(k)} Q_{|X_0|}^2(u) du < \infty \, ,\eeq 
where $Q_{|X_0|}$ is given in Definition \ref{defquant}. Then the series
$\bkE (X_0^2) + 2 \sum_{k \geq 1} \bkE (X_0 X_k)$ is convergent to some nonnegative real $\sigma^2$ and one can construct a sequence $(Z_i)_{i \geq 1}$ of zero mean i.i.d. Gaussian variables with variance $\sigma^2$ such that (\ref{strassen}) holds true with $a_n = ( n \log \log n)^{1/2}$. As shown in Theorem 3 of Rio (1995-a), the condition (\ref{quantrio}) cannot be improved. Recently Dedecker, Gou\"{e}zel and Merlev\`{e}de (2010) proved that this result still holds if we replace the Rosenblatt strong mixing coefficients  $\alpha(n)$ by the weaker coefficients defined in (\ref{defalpha}), provided that the underlying sequence is ergodic.

Still in the strong mixing setting, the best extension, up to our knowledge, of the Koml\'os, Major and Tusn\'ady results  is due to Shao and Lu (1987). Applying the Skorohod embedding, they obtained the following result (see also 
Corollary 9.3.1 in Lin and Lu (1996)): Let $p \in ]2,4]$ and $r >p$. Assume that
\beq \label{Shaocond}
{\mathbb E}(|X_0|^{r})< \infty \quad \text{ and} \quad \sum_{n \geq 1} (\alpha(n))^{(r-p)/(rp)} < \infty \, .
\eeq 
Then the series
$\bkE (X_0^2) + 2 \sum_{k \geq 1} \bkE (X_0 X_k)$ is convergent to some nonnegative real $\sigma^2$ and  one can construct a sequence $(Z_i)_{i \geq 1}$ of zero mean i.i.d. Gaussian variables with variance $\sigma^2$ such that (\ref{strassen}) holds true with $a_n =n^{1/p}(\log  n)^{1 + (1 + \lambda)/p}$, where $\lambda = (\log 2)/\log (r / (r-2))$.

\medskip

Comparing (\ref{Shaocond}) with (\ref{quantrio}) when $p$ is close to $2$, there appears to be a gap between the two above results.  A reasonable conjecture is that 
Shao and Lu's result  still holds under the weaker condition
\beq \label{ros}
{\mathbb E}(|X_0|^{p})< \infty \quad \text{ and} \quad \sum_{k=1}^\infty k^{p-2}\int_0^{\alpha(k)} Q_{|X_0|}^p(u) du < \infty \, ,
\eeq 
since the Rosenthal inequality of order $p$ is true under (\ref{ros}) (see Theorem 6.3 in Rio (2000)) and may fail to hold if this condition is not satisfied (see Rio (2000), chapter 9).  To compare (\ref{ros}) with (\ref{Shaocond}), note that (\ref{ros}) is implied by: for $r>p$,
$$
\sup_{x >0} x^r {\mathbb P}(|X_0| > x) < \infty \quad \text{ and } \quad \sum_{n=1}^\infty 
n^{p-2} (\alpha(n) )^{(r-p)/r} < \infty\, ,
$$
which is much weaker than (\ref{Shaocond}). For example, in the case of bounded random variables ($r = \infty$),  (\ref{Shaocond}) needs $\alpha (n) = O (n^{-p})$, 
while (\ref{ros}) holds as soon as $\alpha (n) = O (n^{1-p} (\log n)^{-1-\eps})$ for some positive $\eps$. 

Let us now give an outline of our results and methods of proofs. Our main result is Theorem \ref{ThSM}, which ensures in particular that, for $p \in ]2,3[$, (\ref{strassen}) holds for $a_n = n^{1/p} (\log n)^{1/2-1/p}$ under (\ref{ros}). Furthermore the error in ${\mathbb L}^2$ is of the same order. The proof of our Theorem \ref{ThSM} is based on an explicit construction of the approximating sequence of i.i.d. Gaussian random variables with the help of conditional quantile transformations. From our construction, the ${\mathbb L}^2$ approximating error 
between dyadic blocks of the initial sequence and the gaussian one can be handled with the help of a conditional version of a functional inequality due to Rio (1998), linking the Wasserstein distance $W_2$ with the Zolotarev distance $\zeta_2$ (see our Proposition \ref{condversion}). 
This method allows us to get a smaller logarithmic factor than the extra factor $(\log n)^{1/2}$ induced by the Skorohod embedding. 
Moreover, it is possible to adapt it (by conditioning up to the future rather than to the past) to deal with the partial sums of non necessarily bounded functions $f$ of iterates of expanding maps such as those considered in Section \ref{dynsys}. For such maps, Theorem \ref{ASmap} completes results obtained by Melbourne and Nicol (2005, 2009) when $f$ is H\"{o}lder continuous. The rest of the paper is organized as follows: Section \ref{proofs} is devoted to the proof of the main results whereas the technical tools are stated and proven in Appendix.

\section{Definitions and main result}\label{MR} Let $(\Omega,{\cal A}, {\mathbb P} )$ be a
probability space.  Assume that there exists some strictly stationary sequence 
$(Y_i)_{i \in {\mathbb Z}}$ of real valued random variables on this probability space, and that the probability space
$(\Omega,{\cal A}, {\mathbb P} )$   is large enough  to contain a sequence  
$(\delta_i)_{i \in {\mathbb Z}}$ of independent random variables with uniform distribution
over $[0,1]$, independent of $(Y_i)_{i \in {\mathbb Z}}$. Define the nondecreasing
filtration $({\cal F}_i)_{i \in {\mathbb Z}}$ by 
${\cal F}_i = \sigma ( (Y_k, \delta_k) : k \leq i)$. Let ${\cal {F}}_{-\infty} = \bigcap_{i \in {\mathbb Z}} {\cal {F}}_{i}$ and ${\cal {F}}_{\infty} = \bigvee_{i \in
{\mathbb Z}} {\cal {F}}_{i}$. We shall denote sometimes by ${\mathbb
E}_i$ the conditional expectation with respect to ${\mathcal F}_i$.

\setcounter{equation}{0} In this section we give rates of
convergence in the almost sure and ${\mathbb L}^2$ invariance principle for functions of a stationary
sequence $(Y_i)_{i \in {\mathbb Z}}$ satisfying weak dependence conditions that we specify below.

\begin{definition} \label{defquant}
 For any nonnegative random variable
$X$, define the ``upper tail'' quantile function $Q_X $ by $
Q_X (u) = \inf \left \{  t \geq 0 : \p \left(X >t \right) \leq
u\right \} $.
\end{definition}
This function is defined on $[0,1]$, non-increasing, right
continuous, and has the same distribution as $X$. This makes it
very convenient to express the tail properties of $X$ using
$Q_X$. For instance, for $0<\varepsilon<1$, if the distribution
of $X$ has no atom at $Q_X(\varepsilon)$, then
  \begin{equation*}
  \bkE ( X \I_{X > Q_X(\varepsilon)})=\sup_{\p(A)\leq \varepsilon} \bkE(X \I_A) = \int_0^\varepsilon Q_X(u) du\, .
  \end{equation*}

\begin{definition}
\label{defclosedenv} Let $\mu$ be the probability distribution of a
random variable $X$. If $Q$ is an integrable quantile function,
let $\tMon( Q, \mu)$ be the set of functions $g$ which are
monotonic on some open interval of ${\mathbb R}$ and null
elsewhere and such that $Q_{|g(X)|} \leq Q$. Let $\tMonm( Q,
\mu)$ be the closure in ${\mathbb L}^1(\mu)$ of the set of
functions which can be written as $\sum_{\ell=1}^{L} a_\ell
f_\ell$, where $\sum_{\ell=1}^{L} |a_\ell| \leq 1$ and $f_\ell$
belongs to $\tMon( Q, \mu)$.
\end{definition}

\begin{definition}
For any integrable random variable $X$, let us write
$X^{(0)}=X- \bkE(X)$.
For any random variable $Y=(Y_1, \cdots, Y_k)$ with values in
${\mathbb R}^k$ and any $\sigma$-algebra $\F$, let
\[
\alpha(\F, Y)= \sup_{(x_1, \ldots , x_k) \in {\mathbb R}^k}
\left \| \bkE \Big(\prod_{j=1}^k (\I_{Y_j \leq x_j})^{(0)} \Big | \F \Big)^{(0)} \right\|_1.
\]
For the sequence ${\bf Y}=(Y_i)_{i \in {\mathbb Z}}$, let \begin{equation}
\label{defalpha} \alpha_{k, {\bf Y}}(0) =1 \text{ and }\alpha_{k, {\bf Y}}(n) = \max_{1 \leq l \leq
k} \ \sup_{ n\leq i_1\leq \ldots \leq i_l} \alpha(\F_0,
(Y_{i_1}, \ldots, Y_{i_l})) \text{ for $n>0$}.
\end{equation}
\end{definition}

For any positive $n$, $\alpha_{k, {\bf Y}}(n) \leq \alpha (n)$, where $\alpha (n)$ is defined by
 (\ref{defalpharosen}). We now introduce some quantities involving the rate of
mixing and the quantile function $Q$. Define  
\begin{equation}
\label{definitionR} 
\alpha^{-1}_{2,{\bf Y}}(x) = \min\{q\in \BBN \tq \alpha_{2,{\bf Y}}(q)\leq x\}
\ \text{ and }\ 
  R(x)=\alpha^{-1}_{2,{\bf Y}}(x)( Q(x) \vee 1)
\end{equation}
(note that $\alpha^{-1}_{2,{\bf Y}}(x) \geq 1$ for $x < 1$). Set, for $p\geq 1$, 
 \begin{equation}
\label{defmoments} 
M_{p,\alpha} (Q) =  \int_0^1 R^{p-1}(u) Q(u) du \ \text{ and }\ 
\Lambda_{p,\alpha} (Q) = \sup_{u\in]0,1]} u R^{p-1} (u) Q(u)  .
\end{equation}
Note that, if $M_{p,\alpha} (Q) <\infty$ then $\Lambda_{p,\alpha} (Q) < \infty$, 
Also, if $\Lambda_{p,\alpha} (Q) < \infty$, then $M_{r,\alpha} (Q) <\infty$ for 
any $r<p$. Let us now state our main result.

\begin{Theorem} \label{ThSM}   Let  $X_i = f(Y_i) - \bkE ( f(Y_i))$ where $f$ belongs to $\tMonm( Q, P_{Y_0})$ (here $P_{Y_0}$ denotes the law of $Y_0$).  
Assume that $M_{2,\alpha} (Q) < \infty$. Then the series
$\bkE (X_0^2) + 2 \sum_{k \geq 1} \bkE (X_0 X_k)$ is convergent to some nonnegative real  $\sigma^2$.
Now let $p \in ]2,3]$ and suppose that  $\Lambda_{p,\alpha} (Q) < \infty$ in the case $p<3$ or $M_{3,\alpha} (Q) < \infty$ in the case $p=3$. 
\begin{enumerate} 
\item 
Assume that $\sigma^2 >0$. Then: 
\begin{enumerate}
\item 
there exists a sequence $(Z_i)_{i \geq 1}$ of iid random variables with law $N(0,\sigma^2)$
such that, setting $\Delta_k = \sum_{i=1}^k (X_i - Z_i)$, 
$$
 \sup_{k \leq n} |\Delta_k| =  O ( n^{1/p} (\log n)^{1/2-1/p}) \text{ in }{\mathbb  L}^2 
\text{ and a.s. for $p< 3$ if } M_{p,\alpha} (Q) < \infty . 
$$
\item For any $\varepsilon >0$, there exists a sequence 
$(\widetilde Z_i)_{i\geq 1}$ of iid  random variables with law $N(0,\sigma^2)$ such that, 
setting  $\widetilde \Delta_k = \sum_{i=1}^k (X_i -\widetilde Z_i)$, 
$$
 \sup_{k \leq n} |\widetilde \Delta_k| =  
 O ( n^{1/p} (\log n)^{1/2}(\log \log n)^{(1+ \varepsilon )/p}) \text{ a.s.}
$$
\end{enumerate}
\item Assume that  $\sigma^2 =0$. Let $S_k = \sum_{i=1}^k X_i$. 
Then
\begin{enumerate}
\item $\sup_{k \leq n} | S_k | =  O ( n^{1/p} )
\text{ in }{\mathbb  L}^2 \text{ and }  \sup_{k \leq n} | S_k | =  O ( (n\log n)^{1/p} 
(\log\log n)^{(1+\varepsilon)/p} )  \text{ a.s.}$
\item If $p<3$ and $M_{p,\alpha} (Q) < \infty$, then 
$\sup_{k \leq n} | S_k | =   o ( n^{1/p} )$ a.s. 
\end{enumerate} 
\end{enumerate} 
\end{Theorem}
\begin{Remark} \label{remarkequicond} The condition $M_{p,\alpha} (Q) < \infty$ can be rewritten in a complete equivalent way as
\beq \label{Condstrong}
 \sum_{k\geq 0} (1\vee k)^{p-2} \int_0^{\alpha_{2, {\bf Y}}(k)} Q^p(u) du
 < \infty \, .\eeq 
(see Annexe C in Rio (2000)), which corresponds to (\ref{ros}) with $\alpha_{2,{\bf Y}} (k)$ 
instead of $\alpha (k)$. 
\end{Remark}
\noindent
{\sl Applications to geometric or arithmetic rates of mixing. } Below we denote by $H$ the cadlag 
inverse of the function $Q$. Assume first that, for some $a$ in $]0,1[$, 
$\alpha_{2 , {\bf Y}} (n) = O ( a^n)$ as $n \rightarrow \infty$.
Then $\alpha_{2 , {\bf Y}}^{-1} (u) = O ( |\log u| )$ as $u$ decreases to $0$. 
Consequently $M_{p,\alpha} ( Q )  < \infty$ as soon as 
$$
\int_0^1 |\log u|^{p-1} Q^p (u) du < \infty . 
$$
This condition holds if
$H(x) = O (\, (x\log x)^{-p} (\log\log x)^{-(1+\eps)} )$ as  $x \rightarrow \infty$.
In  a similar way $\Lambda_{p,\alpha} (Q) < \infty$ if one of the following equivalent weaker conditions holds:
\begin{eqnarray*}
Q (u) = O ( u^{-1/p} |\log u|^{-1+ (1/p)} ) \hbox{ as } u\downarrow 0 \, ,
\
H(x) = O (\, x^{-p} (\log x)^{1-p} ) \hbox{ as } x\uparrow \infty .
\end{eqnarray*}
\par
Suppose now that, for some real $q>2$,  $\alpha_{2 , {\bf Y}} (n) = O ( n^{1-q} )$ as $n \rightarrow \infty$.
Then $\alpha_{2 , {\bf Y}}^{-1} (u) = O ( u^{-1/(q-1)} )$ as $u\rightarrow 0$. For $p$ in $[2,q[$,   
we get that $M_{p,\alpha} ( Q )  < \infty$ as soon as 
$$
\int_0^1 |u|^{-1/(q-1)} Q^p (u) du < \infty . 
$$
This condition holds if $H(x) = O (\, (x^p \log (x) (\log\log x)^{1+\eps} )^{-(q-1)/(q-p)} )$ as $x\rightarrow \infty$.
In  a similar way $\Lambda_{p,\alpha} (Q) < \infty$ if and only if
$H(x) = O (\, x^{-p(q-1)/(q-p)} )$ as $x \rightarrow \infty$. 
Note also that $\Lambda_{q,\alpha} (Q) < \infty$ if and only if $Q$ is uniformly bounded over $]0,1]$. 
\par

\section{Application to dynamical systems} \label{dynsys}
\setcounter{equation}{0}
In this section, we consider a class of piecewise expanding
maps $T$ of $[0,1]$ with a neutral fixed point, and their
associated Markov chain $Y_i$ whose transition kernel is the
Perron-Frobenius operator of $T$ with respect to the absolutely
continuous invariant probability measure. Applying Theorem \ref{ThSM}, we give a large class
of unbounded functions $f$ for which we can give rates of convergence close to optimal in the strong invariance principle of  the partial sums of
both $f\circ T^i$ and $f(Y_i)$. 

For  $\gamma$ in $]0, 1[$, we consider the intermittent map
$T_\gamma$ from $[0, 1]$ to $[0, 1]$, which is a modification of the
Pomeau-Manneville map (1980):
$$
   T_\gamma(x)=
  \begin{cases}
  x(1+ 2^\gamma x^\gamma) \quad  \text{ if $x \in [0, 1/2[$}\\
  2x-1 \quad \quad \quad \ \  \text{if $x \in [1/2, 1]$} \, .
  \end{cases}
$$
We denote by $\nu_\gamma$ the unique $T_\gamma$-invariant
probability measure on $[0, 1]$ which is absolutely continuous with
respect to the Lebesgue measure. We denote by $K_\gamma$ the
Perron-Frobenius operator of $T_\gamma$ with respect to
$\nu_\gamma$. Recall that for any bounded measurable functions $f$ and $ g$,
$$
\nu_\gamma(f \cdot  g\circ T_\gamma)=\nu_\gamma(K_\gamma(f) g) \, .
$$
Let $(Y_i)_{i \geq 0}$ be a stationary Markov chain with invariant
measure $\nu_\gamma$ and transition Kernel $K_\gamma$. It is well
known (see for instance Lemma XI.3 in Hennion and Herv\'e (2001))
that on the probability space $([0, 1], \nu_\gamma)$, the random
variable $(T_\gamma, T^2_\gamma, \ldots , T^n_\gamma)$ is
distributed as $(Y_n,Y_{n-1}, \ldots, Y_1)$.

To state our results for those intermittent  maps, we need preliminary definitions. 
\begin{Definition}
A function $H$ from ${\mathbb R}_+$ to $[0, 1]$ is a tail
function if it is non-increasing, right continuous, converges
to zero at infinity, and $x\rightarrow x H(x)$ is integrable.
\end{Definition}
\begin{Definition}
\label{defMonGPM}
If $\mu$ is a probability measure on $\mathbb R$ and $H$ is a
tail function, let $\Mon(H, \mu)$ denote the set of functions
$f:{\mathbb R}\to {\mathbb R}$ which are monotonic on some open interval and null
elsewhere and such that $\mu(|f|>t)\leq H(t)$. Let $\Monm(H,
\mu)$ be the closure in ${\mathbb L}^1(\mu)$ of the set of
functions which can be written as $\sum_{\ell=1}^L a_\ell
f_\ell$, where $\sum_{\ell=1}^L |a_\ell| \leq 1$ and $f_\ell\in
\Mon(H, \mu)$.
\end{Definition}

Note that a function belonging to $\Monm(H, \mu)$ is allowed to
blow up at an infinite number of points. Note also that any
function $f$ with bounded variation ($\BV$) such that $|f|\leq
M_1$ and $\|df\|\leq M_2$ belongs to the class $\Monm(H, \mu)$
for any $\mu$ and the tail function $H=\I_{[0, M_1+2M_2)}$
(here and henceforth, $\|df\|$ denotes the variation norm of
the signed measure $df$). 
In the unbounded case, if a function $f$ is piecewise monotonic with $N$ branches, then 
it belongs to $\Monm(H, \mu)$ for $H(t)=\mu(|f|>t/N)$. Finally, let us
emphasize that there is no requirement on the modulus of
continuity for functions in $\Monm(H, \mu)$.

Let $Q$ denote the cadlag inverse of $H$. Then, for the random variable $X$ defined by $X(\omega) = 
\omega$, $\Mon(H, \mu) = \tMon( Q, \mu)$ and $\Monm(H, \mu) = \tMonm ( Q, \mu)$. Furthermore
Proposition 1.17 in Dedecker, Gou\"{e}zel and 
Merlev\`{e}de (2010) states that there exists a positive constant $C$ such 
that, for any $n>0$, $\alpha_{2,{\bf Y}}(n) \leq C n^{(\gamma -1)/\gamma}$. In addition,
the computations page 817 in the same paper show that, for $p\gamma<1$,  the  integrability conditions below are equivalent:
\begin{equation} \label{ratecond}
 \int_0^1 R^{p-1}(u) Q(u) du < \infty \ \hbox{ and }\  \int_0^{\infty} x^{p-1} (H(x))^{\frac{1-p\gamma}{1-\gamma}} dx < \infty \, .
\end{equation}
 Also, for $p$ in $]2, 1/\gamma[$, 
\begin{equation} \label{ratecond2}
\Lambda_{p,\alpha} (Q) < \infty \ \hbox{if and only if}\ H(x) = O (x^{-p(1-\gamma)/(1-p\gamma)} )
\ \hbox{as}\ x\rightarrow\infty \, 
\end{equation}
and, for  $p= 1/\gamma$ and $H =\I_{[0, M)}$, $\Lambda_{p,\alpha} (Q) <\infty$ (see the previous section). 

A modification of the proof of Theorem \ref{ThSM}  leads to the result below for the  
 Markov chain or the dynamical system associated to the transformation $T_\gamma$.

\begin{Theorem} \label{ASmap}  Let $\gamma <1/2$.
 Let $f \in {\mathcal F}(H, \nu_\gamma)$ for some tail function $H$  satisfying (\ref{ratecond})
with $p=2$. Then the series 
\begin{equation}\label{seriecov}
\sigma^2(f)= \nu_\gamma((f-\nu_\gamma(f))^2)+ 2
\sum_{k>0} \nu_\gamma ((f-\nu_\gamma(f))f\circ T_\gamma^k)
\end{equation}
converges absolutely to some nonnegative number $\sigma^2 (f)$. Let $p \in ]2,3]$  satisfying 
$p\leq 1/\gamma$. Let $Q$ denote the cadlag inverse of $H$. Suppose that  
$\Lambda_{p,\alpha} (Q) < \infty$ in the case $p<3$ or $M_{3,\alpha} (Q) < \infty$ in the case $p=3$. 
\begin{enumerate}
\item Let $(Y_i)_{i \geq 1}$ be a stationary Markov chain with
transition kernel $K_{\gamma}$ and invariant measure $\nu_\gamma$, and let
$X_i=f(Y_i)-\nu_\gamma(f)$. The sequence $(X_i)_{i \geq 0}$ satisfies the conclusions of Items 1 and 2 of Theorem \ref{ThSM} with $\sigma^2 = \sigma^2(f)$. 
\item If $\sigma^2(f) = 0$, the sequence $(f\circ T_\gamma^i- \nu_\gamma(f) )_{i \geq 1}$ satisfies the conclusions of Item 2 of Theorem \ref{ThSM}. If $\sigma^2(f) > 0$, enlarging the probability space $([0, 1], \nu_\gamma)$, there exist  sequences $(Z^*_i)_{i \geq
1}$  and $(\tilde Z^*_i)_{i \geq
1}$ of iid random variables with law $N( 0, \sigma^2 (f) )$ such that the random variables 
$\Delta_k = \sum_{i=1}^k (f\circ T_\gamma^i- \nu_\gamma(f) - Z_i^* )$ 
satisfy the conclusions of Item 1(a) of Theorem \ref{ThSM} and the random variables 
$\tilde \Delta_k =\sum_{i=1}^k (f\circ T_\gamma^i- \nu_\gamma(f) - \tilde Z_i^* )$ satisfy the conclusion of Item 1(b).
\end{enumerate}
\end{Theorem}

Item 1 is direct by using Theorem \ref{ThSM} together with (\ref{ratecond}) and (\ref{ratecond2}).  Item 2 requires a proof that is given in Section \ref{proofASmap}.

\begin{Remark}
Theorem \ref{ASmap} can be extended to generalized Pomeau-Manneville
map (or GPM map) of parameter $\gamma \in (0,1)$ as defined in Dedecker, Gou\"{e}zel and 
Merlev\`{e}de (2010).
\end{Remark}

In the specific case of bounded variation functions, Theorem \ref{ASmap} provides the almost sure
invariance principle below for the dynamical system associated to $T_\gamma$. Below we give the results
in the case $\sigma^2 (f)>0$. The rates are slightly better in the case $\sigma^2 (f) = 0$. 

\begin{Corollary} \label{ASmapB}  Let $\gamma \in ]1/3 , 1/2[$ and $f$ be 
a function of bounded variation. Then the series in (\ref{seriecov})
converges absolutely to some nonnegative number $\sigma^2 (f)$ and, for any $\varepsilon >0$, 
there exists a sequence 
$(\widetilde Z_i^*)_{i\geq 1}$ of iid  random variables with law 
$N(0,\sigma^2 (f) )$ such that 
$$
 \sup_{k \leq n} | \sum_{i=1}^k (f\circ T_\gamma^i- \nu_\gamma(f) - \tilde Z_i^* ) | =  
 O ( n^{\gamma} (\log n)^{1/2}(\log \log n)^{(1+ \varepsilon )\gamma}) \text{ a.s.}
$$
\end{Corollary}

For the maps under consideration and  H\"older continuous functions $f$, by using an approximation argument introduced by Berkes and Philipp (1979), 
Melbourne and Nicol (2009) obtained the following explicit error term in the almost sure invariance principle (see their Theorem 1.6 and their Remark 1.7): Let $p >2$ and $0<\gamma< 1/p$, then the error term in the almost sure invariance results is $O(n^{\beta + \varepsilon})$ where $\varepsilon >0$ is arbitrarily small and $\beta = \frac{\gamma}{2} + \frac{1}{4}$ if $\gamma$ belongs to $]1/4, 1/2[$
and $\beta = \frac{3}{8} $ if $\gamma \leq 1/4$. Consequently, for the modification of the
Pomeau-Manneville map and functions $f$ of bounded variation, Corollary \ref{ASmapB} improves 
the error in the almost sure invariance principle obtained in Theorem 1.6 in Melbourne and Nicol (2009). Note also that, for $\gamma < 1/3$ and $f$ 
of bounded variation, condition (\ref{ratecond}) is satisfied with $p=3$,  and Theorem \ref{ASmap}
gives the error term $O(n^{1/3} (\log n)^{1/2} (\log\log n)^{(1+ \varepsilon )/3})$
in the almost sure invariance principle. 

\section{Proofs} \label{proofs}

\setcounter{equation}{0}
From now on, we denote by $C$ a numerical constant  which may vary
from line to line. Throughout the proofs, to shorten the notations, we write 
$\alpha (n) = \alpha_{2, {\bf Y}} (n)$ and $\alpha^{-1} (u) = \alpha_{2, {\bf Y}}^{-1}  (u)$.
We also set, for $\lambda >0$, 
\begin{equation}\label{3alphatronq}
M_{3,\alpha} (Q, \lambda) = \int_0^1 Q(u) R(u) (R(u) \wedge \lambda) du . 
\end{equation}
We start by recalling some fact proved in Rio (1995-b), Lemma A.1.: for $p$ in $]2,3[$, 
\begin{equation} \label{RioA1} 
M_{3,\alpha} (Q, \lambda) = O ( \lambda^{3-p} ) \ \hbox{ as }\  \lambda \rightarrow +\infty 
\ \hbox{ if }\  \Lambda_{p,\alpha} (Q) < \infty . 
\end{equation}

\subsection{Proof of Theorem \ref{ThSM}.}

Assume first that $\sigma^2 >0$. For $L \in {\mathbb N}$, let $m(L) \in {\mathbb N}$ be such that $m(L)\leq L$. Let 
\begin{equation*} \label{defUkL}
I_{k,L} = ]2^L + (k-1)2^{m(L)} ,  2^L + k 2^{m(L)}] \cap {\mathbb N} \ \text{and}\ U_{k,L} = \sum_{i\in I_{k,L}} X_i \, ,
\, k \in \{1, \cdots, 2^{L-m(L)} \} \, .
\end{equation*}
For $k \in \{1, \cdots, 2^{L-m(L)} \} $, let $V_{k,L}$ be the
${\mathcal N}(0, \sigma^2 2^{m(L)} )$-distributed random variable
defined from $U_{k,L}$ via the conditional quantile transformation,
that is \beq \label{defVkN} 
V_{k,L} = \sigma 2^{m(L)/2} \Phi^{-1}
(\widetilde F_{k,L}( U_{k,L}- 0 ) + \delta_{2^L + k 2^{m(L)}} (\widetilde
F_{k,L}( U_{k,L} ) - \widetilde F_{k,L} (U_{k,L} - 0) ) )  \, ,
\eeq where
$\widetilde F_{k,L}:= F_{U_{k,L} | {\mathcal F}_{2^{L}+ (k-1)2^{m(L)}}}$ is the d.f. of 
$P_{U_{k,L} | {\mathcal F}_{2^L + (k-1)2^{m(L)}
}}$ (the conditional law of $U_{k,L}$ given $\mathcal{F}_{2^L +
(k-1)2^{m(L)} }$)     and $\Phi^{-1}$ the inverse of the
standard Gaussian distribution function $\Phi$. Since $\delta_{2^L + k 2^{m(L)}}$ is
independent of  ${\mathcal F}_{2^L + (k-1)2^{m(L)}}$, the random variable $V_{k,L}$
is independent of ${\mathcal F}_{2^L + (k-1)2^{m(L)}}$, and has the Gaussian distribution
$N (0, \sigma^2 2^{m(L)})$. By induction on $k$, the random variables $(V_{k,L})_k$
are mutually independent and independent of ${\mathcal F}_{2^L}$. 
In addition
 \begin{eqnarray*} \label{transquantcond}\bkE
(U_{k,L} - V_{k,L})^2 & = & \bkE \int_0^1 \big ( F^{-1}_{U_{k,L} |
{\mathcal F}_{2^L + (k-1)2^{m(L)} }}(u)- \sigma 2^{m(L)/2}\Phi^{-1}(u)
\big )^2 du \\ & : = &  \bkE  \big ( W_2^2 (P_{U_{k,L} | {\mathcal
F}_{2^L + (k-1)2^{m(L)} }} , G_{\sigma^2 2^{m(L)} }) \big ) \, ,
\end{eqnarray*} 
where $G_{\sigma^2 2^{m(L)} }$ is the Gaussian distribution
$N (0, \sigma^2 2^{m(L)})$. Using Proposition
\ref{condversion} and stationarity, we then get that there exists a positive constant
$C$ such that 
\beq\label{majw2cond} \bkE (U_{k,L} - V_{k,L})^2 \leq
C 2^{m(L)/2} M_{3,\alpha} (Q ,2^{m(L)/2})  \, . 
\eeq
\par
Now we construct a
sequence $(Z'_i)_{i \geq 1}$ of i.i.d. Gaussian random variables with
zero mean and variance $\sigma^2$  as follows.  Let $Z'_1 = \sigma \Phi^{-1} (\delta_1)$. For any
$L \in {\mathbb N}$ and any $k \in \{1, \cdots, 2^{L -m(L)} \}$ 
the random variables $(Z'_{2^L + (k-1)2^{m(L)} +1}, \ldots , Z'_{2^L + k 2^{m(L)}})$ 
are defined in the following way. If $m(L)=0$, then 
$Z'_{2^L + k 2^{m(L)}} = V_{k,L}$.
If $m(L)>0$, then by the Skorohod lemma (1976), there exists some measurable function $g$ from ${\mathbb R} \times [0,1]$ 
in ${\mathbb R}^{2^{m(L)}}$ such that, for any pair 
$(V,\delta)$ of independent random variables with respective laws $N(0,\sigma^2 2^{m(L)})$ and 
the uniform distribution over $[0,1]$, $g(V,\delta) = (N_1, \ldots N_{2^{m(L)}})$ is a Gaussian random vector with i.i.d. components such that $V = N_1 + \cdots + N_{2^{m(L)}}$. 
We then set 
$$
(Z'_{2^L + (k-1)2^{m(L)} +1}, \ldots , Z'_{2^L + k 2^{m(L)}})  = g (V_{k,L} , \delta_{2^L + (k-1)2^{m(L)} +1} ). 
$$
The so defined sequence $(Z'_i)$ has the prescribed distribution. 
\par\medskip
Set $S_j = \sum_{i=1}^j X_i$ and $T_j = \sum_{i=1}^j Z'_i$. Let 
\begin{eqnarray*} \label{0dec} 
D_L := \sup_{ \ell \leq 2^{L}} | \sum_{i = 2^L +1}^{2^L + \ell} (X_i -Z'_i)| \, .
\end{eqnarray*}
Let $N\in {\BBN}^*$ and let $k  \in ]1, 2^{N+1}]$.
We first notice that $D_L \geq | (S_{2^{L + 1} } -  T_{2^{L + 1} }) - ( S_{2^L } - T_{2^L } )|$,
so that, if $K$ is the  integer such that $2^K < k \leq  2^{K+1}$, $|S_k - T_k| \leq |X_1 - Z'_1| + D_0 + D_1 + \cdots + D_K$. 
Consequently since $K \leq N$,
\begin{eqnarray} \label{1dec} \sup_{1 \leq k \leq 2^{N+1} }|S_k -T_k|  \leq 
|X_1 - Z'_1| + D_0 + D_1 + \cdots + D_N . 
\end{eqnarray}
\par
We first notice that the following decomposition is valid: 
\beq \label{decsup} 
D_L  \leq D_{L,1} + D_{L,2} \, ,
 \eeq 
where 
$$ 
D_{L,1}:= \sup_{k \leq 2^{L-m(L)} } \Big|
\sum_{ \ell =1}^k (U_{\ell,L} - V_{\ell,L}) \Big|
\ \text{and}\ 
D_{L,2}:= \sup_{ k \leq 2^{L-m(L)}} \sup_{ \ell \in I_{k,L} } 
\Big| \sum_{i = \inf I_{k,L}  }^{\ell} (X_i -Z_i) \Big| .
$$ 
The main tools for proving Theorem \ref{ThSM} will be the two lemmas below.
The first lemma allows us to control the fluctuation term $D_{L,2}$. 

\begin{Lemma} \label{LemDL2} There exists positive constants $c_1$, $c_2 \geq 2$, $c_3$ and $c_4$ such that, 
for any positive $\lambda$,  
\begin{equation}
\BBP  ( D_{L,2} \geq 2\lambda )  \leq  (c_1+2) 2^L \exp \Bigl( - \frac{\lambda^2}{c_2 \sigma^2
2^{m(L)}} \Bigr) + 2^L\lambda^{-3} \bigl(  c_3 
M_{3,\alpha} (Q,\lambda) + c_4 \sigma^3 \bigr)  \, . \label{apppropineg}
\end{equation}
\end{Lemma}

The second lemma gives a bound in ${\mathbb L}^2$ on the Gaussian approximation term $D_{L,1}$. 

\begin{Lemma} \label{LemDL1} Let $p\in ]2,3]$. Suppose that $\Lambda_{p,\alpha} (Q) < \infty$
in the case $p<3$ and $M_{3,\alpha} (Q) < \infty$ in the case $p=3$. Then 
\begin{equation} \label{DL1borneL2}
\|D_{L,1}\|^2_2 \leq C 2^L \Bigl ( 2^{(2-p)m(L)} + 2^{-m(L)/2} M_{3,\alpha} (Q , 2^{m(L)/2})  \Bigr) \, .
\end{equation}
\end{Lemma}

\noindent
{\bf Proof of Lemma \ref{LemDL2}.}  By the triangle inequality together with the stationarity of the sequences $(X_i)_i$ and $(Z_i)_i$,  for any positive $\lambda$,
\beq \label{decdl2} {\mathbb P}
( D_{L,2} \geq 2 \lambda ) \leq 2^{L-m(L)} {\mathbb P} \Bigl(  \sup_{
\ell \leq 2^{m(L)}} | S_{\ell} | \geq \lambda \Bigr) +  2^{L-m(L)} {\mathbb P} 
\Bigl(  \sup_{\ell \leq 2^{m(L)}} | T_{\ell} | \geq \lambda
\Bigr)  \, .\eeq
By L\'evy's inequality (see for instance Proposition 2.3 in Ledoux and Talagrand (1991)),
\begin{equation} {\mathbb P} \Bigl(  \sup_{
 \ell\leq 2^{m(L)}} | T_{\ell} | \geq \lambda
\Bigr) \leq  2 \exp  \Bigl( - \frac{\lambda^2}{2 \sigma^2 2^{ m(L)}
} \Bigr)   \, . \label{LI}
\end{equation}
On the other hand, applying Proposition \ref{inegamax}, we get that
$$
\BBP  \Bigl(  \sup_{ \ell \leq 2^{m(L)}} | S_{\ell} | \geq \lambda \Bigr)  \leq  
c_1\exp \Bigl( - \frac{\lambda^2}{c_2 \sigma^2
2^{m(L)}} \Bigr) + 2^{m(L)}\lambda^{-3} \bigl(  c_3 
M_{3,\alpha} (Q,\lambda) + c_4 \sigma^3 \bigr)   \, .
$$
Collecting the above inequalities, we then get Lemma \ref{LemDL2}. 
\par\medskip\noindent
{\bf Proof of Lemma \ref{LemDL1}.} For any $\ell \in \{ 1, \cdots, 2^{L-m(L)} \} $, let
$\widetilde U_{\ell,L} =U_{\ell,L} - \bkE_{2^L + (\ell-1)2^{m(L)}} (U_{\ell,L})$. 
Then $(\widetilde U_{\ell,L})_{\ell \geq 1}$ is a strictly
stationary sequence of martingale differences adapted to the
filtration $({\mathcal F}_{2^L + \ell 2^{m(L)}})_{\ell \geq 1}$. Notice first that
\beq \label{dec0item2}
\|D_{L,1}\|_2 \leq \Vert \sup_{k \leq 2^{L-m(L)}} | \sum_{ \ell =1}^k (\widetilde
U_{\ell,L} - V_{\ell,L})  | \Vert_2 +\Vert \sup_{k \leq 2^{L-m(L)}}  | \sum_{ \ell =1}^k (\widetilde
U_{\ell,L} - U_{\ell,L})   | \Vert_2  \, .
\eeq
Let us deal with the first term on right hand. Since $V_{\ell,L}$ is independent of 
${\mathcal F}_{2^L + (\ell - 1) 2^{m(L)} }$, the sequence $(\widetilde
U_{\ell,L} - V_{\ell,L})_\ell$ is a martingale difference sequence with respect to the
nondecreasing filtration $({\mathcal F}_{2^L + \ell 2^{m(L)}})_\ell$. Hence, by the 
Doob-Kolmogorov maximal inequality, we get that 
\begin{eqnarray*}
\Vert \sup_{k \leq 2^{L-m(L)}} \big | \sum_{ \ell =1}^k (\widetilde
U_{\ell,L} - V_{\ell,L}) \big  | \Vert^2_2 & \leq & 4
\sum_{\ell=1}^{2^{L-m(L)}} \Vert \tilde  U_{\ell, L} -  V_{\ell, L}
\Vert^2_2 \\
& \leq & 8 \sum_{\ell=1}^{2^{L-m(L)}}  \Vert \widetilde
U_{\ell,L} -  U_{\ell,L} \Vert^2_2 +8 \sum_{\ell=1}^{2^{L-m(L)}}
\Vert U_{\ell, L} - V_{\ell, L}\Vert^2_2 \, .
\end{eqnarray*}
Since $V_{\ell,N}$ is independent of  $\mathcal{F}_{2^L +
(\ell-1)2^{m(L)} }$, $\bkE_{2^L + (\ell-1)2^{m(L)}} (V_{\ell,L}) =0$.
Consequently,
$$ \Vert \tilde  U_{\ell, L} -  U_{\ell, L} \Vert^2_2 = \Vert \bkE_{2^L + (\ell-1)2^{m(L)}} (U_{\ell, L} - V_{\ell,L}) \Vert_2^2
\leq \Vert U_{\ell, L} - V_{\ell,L}\Vert_2^2 \, . $$ Using (\ref{majw2cond}), it follows that
\beq \label{maj5DR}
\Vert \sup_{k \leq 2^{L-m(L)}} | \sum_{ \ell =1}^k (\widetilde
U_{\ell,L} - V_{\ell,L})  | \Vert_2^2 \leq C 2^{L-m(L)/2} M_{3,\alpha} (Q, 2^{m(L)/2} )\, .
\eeq
We  deal now  with the second term in the right hand side of (\ref{dec0item2}). According to Dedecker and Rio's maximal inequality (2000, Proposition 1), we obtain that
\begin{eqnarray} \label{inegadedrio}
& & \Vert \sup_{k \leq 2^{L-m(L)}}  | \sum_{ \ell =1}^k (\widetilde
U_{\ell,L} - U_{\ell,L})   | \Vert_2^2   \leq  4 \sum_{k=1}^{2^{L-m(L)}} \Vert\bkE_{2^L + (k-1)2^{m(L)}} (U_{k,L}) \Vert_2^2
\nonumber \\
& & \quad \quad + 8 \sum_{k=1}^{2^{L-m(L)} -1 }
 \Vert\bkE_{2^L + (k-1)2^{m(L)}} (U_{k,L}) \big ( \sum_{i= k +1}^{2^{L-m(L)}}\bkE_{2^L + (k-1)2^{m(L)}} (U_{i,L}) \big ) \Vert_1 \, .
\end{eqnarray}
Stationarity leads to
\beq \label{maj1DR}
\Vert\bkE_{2^L + (k-1)2^{m(L)}} (U_{k,L}) \Vert_2^2= \Vert\bkE_{0} (S_{2^{m(L)}}) \Vert_2^2 \leq 2 \sum_{i=1}^{2^{m(L)}} \sum_{j=1}^i {\mathbb E} | X_j {\mathbb E}_0 (X_i)| \, .
\eeq
Using Lemma 4 (page 679) in Merlev\`ede and Peligrad (2006), we get that 
$$
{\mathbb E} | X_j {\mathbb E}_0 (X_i)| \leq 3 \int_{0}^{\Vert{\mathbb E}_0 (X_i) \Vert_1} Q_{|X_0|} \circ G_{|X_0|} (u) du \, ,
$$
where $G_{|X_0|}$ is the inverse of $L_{|X_0|}(x) =
\int_0^{x} Q_{|X_0|} (u) du $. We will denote by $L$ and $G$
the same functions constructed from $Q$. Assume first that $X_i = f(Y_i)-\mathbb E(f(Y_i))$ with
$f=\sum_{\ell=1}^L a_{\ell} f_{\ell}$, where $f_\ell \in
\tMon( Q,P_{Y_0})$ and $\sum_{\ell=1}^L |a_\ell| \leq 1$. According
to Proposition \ref{covineg},
  \begin{equation}
  \label{majgamma2}\Vert{\mathbb E}_0 (X_i) \Vert_1 \leq 8 \int_0^{\alpha (i)} Q(u) du \, .
  \end{equation}
Since $Q_{|X_0|}(u)\leq Q_{|f(Y_0)|} (u) + |\mathbb E(f(Y_0))|$, we see that
$\int_0^x Q_{|X_0|}(u) du \leq 2 \int_0^x Q_{|f(Y_0)|}(u) du$.
Since $f=\sum_{\ell=1}^L a_\ell f_\ell$, we get, according to Item (c) of
Lemma 2.1 in Rio (2000),
\[
\int_0^{x} Q_{|X_0|} (u) du \leq 2\sum_{\ell=1}^L \int_0^{x} Q_{|a_\ell
f_\ell(X_0)|} (u) du \leq 2\sum_{\ell=1}^L |a_\ell| \int_0^{x} Q (u) du \, .
\]
Since $\sum_{\ell=1}^L |a_\ell| \leq 1$, it follows that $
G(u/2) \leq G_{|X_0|} (u)$. In particular, $G_{|X_0|}(u)\geq
G(u/8)$. Using the fact that $Q_{|X_0|}$ is non-increasing and the change of variables $w = G(v)$,
  \begin{align*}
 \int_{0}^{\Vert{\mathbb E}_0 (X_i) \Vert_1} Q_{|X_0|} \circ G_{|X_0|} (u) du 
  &\leq \int_0^{\Vert{\mathbb E}_0 (X_i) \Vert_1} Q_{|X_0|} \circ G(u/8) du
  = 8 \int_0^{\Vert{\mathbb E}_0 (X_i) \Vert_1/8} Q_{|X_0|}\circ G(v)dv
  \\&
  =8 \int_0^{ G(\Vert{\mathbb E}_0 (X_i) \Vert_1/8)} Q_{|X_0|}(w) Q(w)dw
  \leq 8 \int_0^{\alpha (i)} Q_{|X_0|} (w) Q(w) dw  \,
  ,
  \end{align*}
where the last inequality follows from (\ref{majgamma2}). Consequently, by Item (c) of
Lemma 2.1 in Rio (2000),
\begin{eqnarray} \label{ineMP}
{\mathbb E} | X_j {\mathbb E}_0 (X_i)| \leq   48 \sum_{\ell=1}^L |a_\ell|  \int_0^{\alpha (i)} Q_{|f_{\ell}(Y_0)|} (u) Q(u) du \leq  48   \int_0^{\alpha (i)} Q^2 (u)du\, ,
    \end{eqnarray}
    and the same inequality holds if
 $f \in \tMonm (Q,P_{Y_0})$ by applying Fatou's lemma. Consequently starting from (\ref{maj1DR}), we derive that
 \beq \label{maj2DR}
\sum_{k=1}^{2^{L-m(L)}} \Vert\bkE_{2^L + (k-1)2^{m(L)}} (U_{k,L}) \Vert_2^2 \leq 96.\, 2^{L-m(L)}\sum_{i=1}^{2^{m(L)}} i \int_0^{\alpha (i)} Q^2 (u)du \, .
\eeq
We now bound up the second term in the right hand side of (\ref{inegadedrio}). Stationarity yields that
$$
 \Vert\bkE_{2^L + (k-1)2^{m(L)}} (U_{k,L}) \big ( \sum_{i= k +1}^{2^{L-m(L)}}\bkE_{2^L + (k-1)2^{m(L)}} (U_{i,L}) \big ) \Vert_1  \leq  \sum_{j=1}^{2^{m(L)}} 
\sum_{i=2^{m(L)}+ 1}^{2^L - (k-1) 2^{m(L)}} {\mathbb E} | X_j {\mathbb E}_0 (X_i)| \, .
$$
Using Inequality (\ref{ineMP}), we then derive that
\begin{equation} \label{maj3DR}
\sum_{k=1}^{2^{L-m(L)} -1 }
 \Vert\bkE_{2^L + (k-1)2^{m(L)}} (U_{k,L}) \big ( \sum_{i= k +1}^{2^{L-m(L)}}\bkE_{2^L + (k-1)2^{m(L)}} (U_{i,L}) \big ) \Vert_1   \leq 48.\, 2^{L} 
\sum_{i=2^{m(L)}+ 1}^{2^L } \int_0^{\alpha (i)} Q^2 (u)du \, .
\end{equation}
Starting from (\ref{inegadedrio}) and considering the bounds (\ref{maj2DR}) and (\ref{maj3DR}), 
we get that
\begin{eqnarray} \label{maj4DR}     \Vert \sup_{k \leq 2^{L-m(L)}}  | \sum_{ \ell =1}^k (\widetilde
U_{\ell,L} - U_{\ell,L})   | \Vert_2^2 & \leq & C 2^{L-m(L)} \!\!  \int_0^1 \!\!Q (u) R (u) 
(\alpha^{-1} (u) \wedge 2^{m(L)} ) du \nonumber \\ 
& \leq & C 2^{L-m(L)} M_{3,\alpha} (Q , 2^{m(L)} ) \, ,
\end{eqnarray}
since $R(u) \geq \alpha^{-1} (u)$.  
Starting from (\ref{dec0item2}) and  considering the bounds (\ref{maj5DR}), 
(\ref{maj4DR}) and  (\ref{RioA1}) in the case $p<3$, we then get (\ref{DL1borneL2}),
which ends the proof of Lemma \ref{LemDL1}. 
\par\medskip

\noindent{\bf Proof of Item 1(a). } We choose $Z_i=Z_i'$ with
\beq \label{defml}
m(L) = \Big [ \frac{2L}{p} -\frac{2}{p} \log_2 L \Big ] \, ,\ \text{so that}\ 
\frac 12 \Bigl(\frac{2^L}{L} \Bigr)^{2/p} \leq 2^{m(L)} \leq \Bigl(\frac{2^L}{L} \Bigr)^{2/p} \, ,
\eeq 
square brackets designating as usual the integer part and $\log_2 (x)= (\log x)/(\log 2)$.  
Starting from (\ref{apppropineg}), we now prove that 
\begin{equation}
D_{L,2} = O ( 2^{L/p} L^{1/2-1/p}) \text{  in } {\mathbb L}^2 
\text{ for $p \leq 3$ and a.s. for $p< 3$ if }\, M_{p,\alpha} (Q) < \infty . \label{DL2}
\end{equation}

To prove the almost sure part in (\ref{DL2}), take
\beq \label{deflambda}\lambda=\lambda_L=  K 2^{m(L)/2} \sqrt{L}
\text{ with } K = \sqrt{2 c_2  \sigma^2 \log 2} \, . \eeq
Then, on one hand,
$$
\sum_{L > 0} 2^L \exp \Bigl( - \frac{\lambda_L^2}{ c_2 \sigma^2 2^{m(L)}} \Bigr) 
= \sum_{L \geq 0} 2^{L-2L} < \infty
\ \text{and}\ 
\sum_{L > 0} 2^{L} \lambda_L^{-3} < \infty \, , 
$$
for $p<3$. On the other hand, since $M_{3,\alpha} (Q, a\lambda) \leq aM_{3,\alpha} (Q,\lambda)$
for any $a\geq 1$, 
\begin{eqnarray*}
 2^L\lambda_L^{-3} M_{3,\alpha} (Q,\lambda_L) \leq
2^{L - 3m(L)/2} L^{-1} M_{3,\alpha} (Q, K 2^{m(L)/2} ) .
\end{eqnarray*}
Consequently, from the choice of $m(L)$ made in (\ref{defml}), 
\begin{equation*}
 \sum_{L > 0} 2^L \lambda_L^{-3}  M_{3,\alpha} (Q,\lambda_L)
\leq C \sum_{L >  0} (2^L / L)^{(p-3)/p}
M_{3,\alpha} (Q, (2^L / L)^{1/p} ).
\end{equation*}
Next, for $p \in ]2,3[$,
$$
 \sum_{ L \, : \, \frac{2^L}{L} \geq R^p(x)} \Big (\frac{2^L}{L}
\Big )^{1 -3/p} \leq C R^{p-3} (x) \text{ and } \sum_{ L \, : \,
\frac{2^L}{L} \leq R^p(x)} \Big (\frac{2^L}{L} \Big )^{1 -2/p} \leq
C R^{p-2} (x) \, ,
$$
which ensures that 
\beq \label{calcul} \sum_{L > 0} 2^L \lambda_L^{-3}  
 M_{3,\alpha} (Q,\lambda_L) \leq C M_{p,\alpha} (Q) . 
\eeq 
Consequently under (\ref{Condstrong}), we derive that $\sum_{L > 0} {\mathbb P} (
D_{L,2} \geq 2 \lambda_L ) < \infty $ implying the almost sure part of (\ref{DL2}) 
via the Borel-Cantelli lemma.
\par\ssk
We now prove the ${\mathbb L}^2$ part of (\ref{DL2}). Clearly
\beq \label{1equalityL2} 
{\mathbb E} ( D_{L,2}^2 ) = 8 \int_0^\infty \lambda {\mathbb P} ( D_{L,2} \geq 2\lambda) d\lambda 
\leq 4 \lambda_L^2 + 8 \int_{\lambda_L}^\infty \lambda {\mathbb P} ( D_{L,2} \geq 2\lambda) d\lambda .
\eeq
We now apply  (\ref{apppropineg}). First, from (\ref{deflambda}), 
$$
\int_{\lambda_L}^\infty \lambda  \exp \Bigl( - \frac{\lambda^2}{c_2 \sigma^2
2^{m(L)}} \Bigr) d\lambda = c_2 \sigma^2 2^{m(L)-L} \ \text{and}\ 
2^L \int_{\lambda_L}^\infty \frac{c_4\sigma^3}{\lambda^2} d\lambda = 
c_4\sigma^3 \frac{2^L}{\lambda_L}.
$$
In the case $p<3$ and $\Lambda_{p,\alpha} (Q) < \infty$, from (\ref{RioA1}),
there exists a positive constant $C$ depending on $p$ and $\Lambda_{p,\alpha} (Q)$ such that 
\beq \label{2inequalityL2} 
 \int_{\lambda_L}^\infty \frac{c_3 2^L}{\lambda^2} M_{3,\alpha} ( Q , \lambda ) 
d\lambda \leq C  \int_{\lambda_L}^\infty  \lambda^{1-p} d\lambda \leq
  \frac{ C 2^L}{(p-2) \lambda_L^{p-2} } \,\, .
\eeq 
Now, by (\ref{deflambda}) again, 
$(K/2) 2^{L/p} L^{1/2 - 1/p} \leq  \lambda_L \leq K 2^{L/p} L^{1/2 - 1/p}$, 
and consequently, collecting the above estimates, we get that 
${\mathbb E} ( D_{L,2}^2 )  = O ( \lambda_L^2 )$, which implies the ${\mathbb L}^2$ part of 
(\ref{DL2}).  
\par\ssk
We now deal with $D_{L,1} $. We will prove that 
\begin{equation}
D_{L,1} = O ( 2^{L/p} L^{1/2-1/p})
\text{  in } {\mathbb L}^2 \text{ for $p \leq 3$ and a.s. for $p< 3$ if }\, 
M_{p,\alpha} (Q) < \infty . \label{DL1}
\end{equation}
  
We first derive from Lemma \ref{LemDL1} that 
$\|D_{L,1}\|^2_2 \leq C 2^{L -m(L) (p-2)/2 }$ (applying (\ref{RioA1}) in the case $p<3$),
which implies the ${\mathbb L}^2$ part of (\ref{DL1}). 
\par\ssk
Next, from (\ref{DL1borneL2}) together with the Markov inequality, 
$$
\sum_{L>0} {\mathbb P} ( D_{L,1} \geq \lambda_L) \leq C \sum_{L>0} 2^{L+ (1-p) m(L)} + 
C \sum_{L>0}  \frac{ 2^L }{L 2^{3m(L)/2} } M_{3,\alpha} (Q , 2^{m(L)/2}) \, ,
$$
where $\lambda_L$ is defined by (\ref{deflambda}). 
Repeating exactly the same arguments as in the proof of (\ref{calcul}), we get that 
the second series on right hand in the above inequality is convergent for $p < 3$. Now 
$2^{L + (1-p) m (L) } \leq 2^{p-1} 2^{L(2-p)/p} L^{2(p-1)/p}$,
which ensures the convergence of the first series on right hand. Hence, by the Borel-Cantelli
lemma $D_{L,1} = O (\lambda_L )$ almost surely, which completes the proof of (\ref{DL1}). 
Finally Item 1(a) of Theorem \ref{ThSM} follows from both (\ref{DL1}), (\ref{DL2}) and (\ref{1dec}) 
and (\ref{decsup}). 
\par\medskip\noindent
{\bf Proof of Item 1(b).}  We choose $\widetilde Z_i= Z_i'$ with 
$m(L) = [ (2L/p) + (2(1 + \varepsilon) / p ) \log_2 (1 \vee \log L)  ]$. 
Following the proof of Item 1(a) with this selection of $m(L)$, Item 1(b) follows.
\par\medskip
\noindent{\bf Proof of Item 2.} Starting from the decomposition (\ref{1dec}), we just have to bound both almost surely and in ${\mathbb L}^2$ the random variables 
$D_L := \sup_{ \ell \leq 2^{L}} | S_{2^L + \ell} - S_{2^L}|$.
Applying Proposition \ref{inegamax} in case where $\sigma^2 =0$, we get that for any positive $\lambda$,
\beq \label{applpropmax0}
\BBP  ( D_L \geq \lambda )  \leq  c 2^{L }\lambda^{-3}  
M_{3,\alpha} (Q , \lambda ) , 
\eeq
where  $c$ is a positive constant. Using computations as in (\ref{1equalityL2}) and (\ref{2inequalityL2}), we then get that for any positive $\lambda_L$,
$\| D_L \|_2^2 \leq 4 \lambda_L^2 + C   2^L \lambda_L^{2-p}$. Choosing $\lambda_L = 2 ^{L/p}$ gives the 
${\mathbb L}^2$ part of Item 2 (a). To prove the almost sure parts, we start from (\ref{applpropmax0}) and choose, for $\delta >0$ arbitrarily small, 
$$
\lambda = 2 ^{L/p} L^{1/p} ( 1 \vee \log L)^{(1+\varepsilon)/p} \text{  and }
\lambda = \delta 2 ^{L/p} \text{ if $p\in ]2,3[$ and } M_{p,\alpha} (Q) < \infty . 
$$
The Borel-Cantelli lemma then implies that almost surely
$$
D_L = O ( 2 ^{L/p} L^{1/p} ( 1 \vee \log L)^{(1+\varepsilon)/p}) 
 \text{ a.s. and } D_L = o (2 ^{L/p})\text{ a.s. if $p\in ]2,3[$ and } M_{p,\alpha} (Q) < \infty \,.
$$
This ends the proof of the almost sure part of Item 2 and then of the theorem.

\subsection{Proof of Item 2 of Theorem \ref{ASmap}.} \label{proofASmap}

If $\sigma^2(f) >0$, similarly as for the proof of Theorem \ref{ThSM}, we start by constructing a sequence $(Z'^*_i)_{i \geq 1}$ of 
i.i.d.~gaussian random variables with mean zero and
variance $\sigma^2(f)$ depending on the sequence $(m(L))_{L \geq 0}$ defined either as in (\ref{defml}) or as in the proof of Item 1(b). 
Define for any $k \in \{1, \cdots, 2^{L-m(L)} \}$, 
$$
I_{k,L} = ]2^L + (k-1) 2^{m(L)} ,   2^L +k 2^{m(L)} ] \cap {\mathbb N} \ \text{and}\ U^*_{k,L} = \sum_{i\in I_{k,L}} (f\circ T_\gamma^i- \nu_\gamma(f))  \, .
$$
For $k \in \{1, \cdots, 2^{L-m(L)} \} $, let $V^*_{k,L}$ be the
${\mathcal N}(0, \sigma^2 2^{m(L)} )$-distributed random variable
defined from $U^*_{k,L}$ via the conditional quantile transformation,
that is \beq \label{defVkNT} 
V^*_{k,L} = \sigma(f) 2^{m(L)/2} \Phi^{-1}
( F^*_{k,L}( U^*_{k,L}- 0 ) + \delta_{2^L +k2^{m(L)} }(
F^*_{k,L}( U^*_{k,L} ) -  F^*_{k,L} (U^*_{k,L} - 0) ) ) \, ,
\eeq where
$ F^*_{k,L}:= F_{U^*_{k,L} | \tilde {\mathcal G}_{ 2^{L}+k2^{m(L)}+1}}$ is the d.f. of 
the conditional law of $U^*_{k,L}$ given $\tilde {\mathcal G}_{2^{L}+k 2^{m(L)}+1}$, 
where $\tilde {\mathcal G}_{m} = \sigma( T_{\gamma}^m, ( \delta_i)_{i \geq m}$ )     and $\Phi^{-1}$ the inverse of the
standard Gaussian distribution function $\Phi$. Since $\delta_{2^L +k2^{m(L)} }$ is
independent of  $\tilde {\mathcal G}_{ 2^L +k2^{m(L)}+1}$, the random variable $V^*_{k,L}$
is independent of $\tilde {\mathcal G}_{2^{L}+k 2^{m(L)}+1}$, and has the Gaussian distribution
$N (0, \sigma^2(f) 2^{m(L)})$. By induction on $k$, the random variables $(V^*_{k,L})_k$
are mutually independent and independent of $\tilde {\mathcal G}_{ 2^{L+1} +1 }$. Let us construct now the
sequence $(Z'^*_i)_{i \geq 1}$  as follows. Let $Z'^*_1= \sigma(f) \Phi^{-1}(\delta_1)$. For any $L \in {\mathbb N}$ and any $k \in \{1, \cdots, 2^{L-m(L)} \}$,
the random variables $(Z'^*_{ 2^L +( k-1) 2^{m(L)}+1}, \ldots , Z'^*_{ 2^L +k2^{m(L)}})$ 
are defined in the following way. If $m(L)=0$, then $
Z'^*_{2^L + k 2^{m(L)}} = V^*_{k,L} $. If $m(L)>0$, then there exists some measurable function $g$ from ${\mathbb R} \times [0,1]$ 
in ${\mathbb R}^{2^{m(L)}}$ such that, for any pair 
$(V,\delta)$ of independent random variables with respective laws $N(0,\sigma^2(f) 2^{m(L)})$ and 
the uniform distribution over $[0,1]$, $g(V,\delta) = (N_1, \ldots N_{2^{m(L)}})$ is a Gaussian random vector with i.i.d. components such that $V = N_1 + \cdots + N_{2^{m(L)}}$. 
We then set 
$$
(Z'^*_{2^L+ ( k-1)2^{m(L)} +1}, \ldots , Z^*_{2^L +k 2^{m(L)}})  = g (V^*_{k,L} ,  \delta_{2^L
+ (k-1) 2^{m(L)}+1} ). 
$$
The so defined sequence $(Z'^*_i)$ has the prescribed distribution. 

\medskip
\par
Set now $ S^*_j= \sum_{i=1}^j (f\circ T_\gamma^i- \nu_\gamma(f))$, $T^*_j = \sum_{i=1}^j Z'^*_i$ if $\sigma^2(f) >0$ and $T^*_j =0$ otherwise, and let 
$$ 
D^*_L := \sup_{ 0 \leq \ell \leq 2^{L}} |(S^*_{2^L + \ell}-T^*_{2^L + \ell}) - ( S^*_{ 2^{L+1} } -T^*_{ 2^{L+1} })| \, .
$$
Similarly as in the proof of (\ref{1dec}), we get that 
\begin{eqnarray} \label{1decdyn} \sup_{1 \leq k \leq 2^{N+1} }|S^*_k -T^*_k|  \leq 
|S^*_1 -T^*_1| + 2 D^*_0 + 2 D^*_1 + \cdots + 2 D^*_N . 
\end{eqnarray}
For any $L \in {\mathbb N}$, on the probability space $([0, 1], \nu_\gamma)$, the random
variable $(T^{2^{L}+1}_\gamma, T^{2^{L}+2}_\gamma, \ldots , T^{2^{L+1}}_\gamma)$ is
distributed as $(Y_{2^{L+1}},Y_{2^{L+1}-1}, \ldots, Y_{2^{L}+1})$, where $(Y_i)_{i \geq 1}$ is a stationary Markov chain with
transition kernel $K_{\gamma}$ and invariant measure $\nu_\gamma$. From our construction of the random variables $Z'^*_i$, for any $L \in {\mathbb N}$, 
$$
(T^{2^{L}+1}_\gamma, \ldots , T^{2^{L+1}}_\gamma,Z'^*_{2^{L}+1}, \ldots , Z'^*_{2^{L+1}}) =^{{\mathcal D}}(Y_{2^{L+1}}, \ldots , Y_{2^{L}+1},Z'_{2^{L+1}}, \ldots , Z'_{2^{L}+1}) \, ,
$$
where the sequence $(Z'_i)_{2^L +1 \leq i \leq 2^{L+1}}$ is defined from $(Y_i,\delta_i)_{2^L < i \leq  2^{L+1}}$ as in the proof of Theorem \ref{ThSM}.
It follows that 
\begin{equation*} \label{equalitylaw*} D^*_L =^{{\mathcal D}} D_L  \text{ where }D_L  := \sup_{ 0 \leq \ell \leq 2^{L}} |(S_{2^L + \ell}-T_{2^L + \ell}) - ( S_{ 2^{L} } -T_{ 2^{L} })| 
\end{equation*}
and, for any $j \geq 1$, $T_j = \sum_{i=1}^j Z'_i$ if $\sigma^2(f) >0$ and $T_j =0$ otherwise. 
Hence we have, for any positive $\lambda$,
$ {\mathbb P} (D^*_{L} \geq \lambda) =  {\mathbb P} (D_{L} \geq \lambda)$.
Proceeding as in the proof of Theorem \ref{ThSM}, Item 2  follows.

\section{Appendix}
\setcounter{equation}{0} 

Next lemma is a parametrized version of Theorem 1 of Rio (1998). We first need the following definition. 
\begin{Definition}
$\Lambda_2$ is the class of
real functions $f$ which are continuously differentiable
and such that $|f'(x) - f'(y) |\leq | x - y |$  for any  $(x,y) \in {\mathbb R} \times {\mathbb R}$.
\end{Definition}

\begin{Lemma} \label{riolma} Let $Z$ be a random variable with values in 
a purely non atomic Lebesgue space $(E, {\mathcal L} (E) , m)$ and ${\mathcal F}  = \sigma (Z)$.
For real random variables $U$ and $V$, let  $P_{U | {\mathcal F}}$ be the law of $U$ 
given ${\mathcal F}$ and $P_V$ be the law of $V$. Assume that $V$ is
independent of $\mathcal F$.  Let $\sigma^2 >0$ and $N$ be a ${\mathcal N} (0,\sigma^2)$-distributed random variable independent of $\sigma( Z,U,V)$. Then
$$
\bkE  \big ( W_2^2 (P_{U | {\mathcal F}} , P_V) \big ) \leq 16 \sup_{f \in \Lambda_2(E)}\bkE \big (   f( U + N , Z) -  f( V + N , Z) \big ) + 8 \sigma^2 
\, ,
$$
where $\Lambda_2 (E)$ denotes the set of measurable functions $f:{\mathbb R}
\times  E \rightarrow {\mathbb R}$ wrt the $\sigma$-fields 
${\mathcal L} ( {\mathbb R} \times E) $ and ${\mathcal B} ({\mathbb R})$, 
such that $f( \cdot, z) \in \Lambda_2$ and $f(0,z)=f'(0,z)=0$  for any 
$z \in E$.
\end{Lemma}

\noindent{\it Proof of Lemma \ref{riolma}.} Notice first that
\beq \label{lissage}
\bkE  \big ( W_2^2 (P_{U | {\mathcal F}} , P_V) \big ) \leq 2 \bkE  \big ( W_2^2 (P_{U + N | {\mathcal F}}  , P_{V + N}) \big ) + 8 \sigma^2 \, .
\eeq Let $G$ be the d.f. of $P_{V + N}$.  
Since $E$ is a Lebesgue space, there exists a regular version of the conditional distribution function of $U + N$ conditionally to $Z$, that is, a function  $(x,  z) \rightarrow F_{z} (x)$ from ${\mathbb R} \times E$ in ${\mathbb R}$ such that, for any real $x$, 
$F_Z (x) =  {\mathbb E} ( \I_{ U + N\leq x} | Z )$ almost surely.

Notice in addition that, for any $z$ in $E$, $F_{z}$ is a $C^{\infty}$ increasing distribution function. 
Let now $H_z (x) = F_z (x) - G(x)$, $ A_z=  \{y \in {\mathbb R} \, : \, H_z (y)= 0 \}$, and for any $(x, z) \in {\mathbb R} \times E$, let
\begin{equation}  \label{defflemma}
h(x, z) = d(x,A_z \cup \{ 0 \} ) \sign H_z (x) \text{ and } f(x ,z ) = \int_0^x h (y , z) dy  \, ,
\end{equation}
where $d(x,A_z \cup \{ 0 \})$ is the distance of $x$ to the random set $A_z\cup \{ 0 \}$ and $\sign y=1$ for $y>0$, $0$ for $y=0$  and $-1$ for $y < 0$. 

For $z$ fixed, $ f(0 ,z )= f'(0 ,z )=0$ and it is shown in Rio (1998, Inequality (7)) that $ f(
\cdot, z)$ belongs to  $\Lambda_2$, and that for any $u \in ]0,1[$, 
$$
  f(F_z^{-1} (u) , z) -  f(G^{-1} (u) , z)   \geq \frac 18 \big (F_z^{-1} (u) - G^{-1} (u) \big )^2 \, ,$$
and therefore that for any $z \in E$, 
\begin{eqnarray} \label{formulemajw2}
W_2^2 (P_{U + N| Z = z } , P_V) & = &  \int_0^1  \big (F_z^{-1} (u) - G^{-1} (u) \big )^2 du \nonumber \\
&  \leq & 8 \Big (  \int_{\mathbb R} f(
x, z) dP_{U+ N | Z=z } - \int_{\mathbb R} f(
x, z) dP_{V+N}  \Big )  \, .
\end{eqnarray}
We prove now that the function $ f$ defined by (\ref{defflemma}) is ${\mathcal L} ( {\mathbb R} \times E) -{\mathcal B} ( {\mathbb R})$ measurable. Notice first that since for any fixed $z$, $x \mapsto h(x, z)$ is continuous we get that 
$$
 f(x ,z ) = \lim_{n \rightarrow \infty} \frac{x}{n} \sum_{i=1}^n h (itn^{-1} , z) \, .
$$
Therefore the mesurability of $ f$ will come from the mesurability of $h$. With this aim, it is 
 enough to prove the mesurability of the restriction $h_n$  of $h$ to $[-n,n] \times E$
for any positive integer $n$. 

Let $\varphi\, : \,  [-n,n] \rightarrow [0,1]$ be the one to one  bicontinuous map defined by 
$\varphi(x) = (n-x)/(2n)$. We then define 
\begin{eqnarray} \label{defggrande}
g \, : \, [0,1]\times E& \rightarrow & {\mathbb R} \nonumber \\
(x, z) & \mapsto & h( \varphi^{-1} (x) ,z ) \, .
\end{eqnarray}
The mesurability of $h_n$ will then follow from the mesurability of $g$.  Since $E$ is purely non atomic, $(E, {\mathcal L} (E), m)$ is isomorph to $([0,1], {\mathcal L} ([0,1]) , \lambda_{[0,1]})$ where ${\mathcal L} ([0,1])$ and $ \lambda_{[0,1]} $ are respectively the Lebesgue $\sigma$-algebra and the Lebesgue measure on $[0,1]$ 
(see for instance Theorem 4.3 in De La Rue (1993)). Consequently the following theorem due to Lipi\'nski (1972) which is recalled in Grande (1976) also holds in $[0,1]\times E$.
\begin{Theorem} (Lipi\'nski (1972)) Let $g$ be a bounded function from $[0,1] \times E$ into ${\mathbb R}$ such that 
\begin{enumerate}
\item the cross sections $g_x (t) = g(x,t)$ and $g^z(t) = g(t,z)$ are respectively ${\mathcal L}(E)$ and  ${\mathcal L}([0,1])$-measurable,
\item for all $t \in [0,1]$, $k_t(z) = \int_0^t g(x,z)dx$ is ${\mathcal L}(E)$-measurable,
\item for all $z \in E$, the cross section $g^z$ is a derivative.
\end{enumerate}
Then $g$ is measurable wrt the $\sigma$-fields  ${\mathcal L} ( [0,1]\times E)$ and ${\mathcal B} ( {\mathbb R})$.
\end{Theorem}
Items 2 and 3 as well as the second part of Item 1 follows directly from the fact that if $z$ is fixed, then the function $x \rightarrow g(x,z)$ is continuous (recall that $h ( \cdot, z)$ and $\varphi^{-1}$ are continuous). It remains to show that for all $x \in [0,1]$ the cross section $g_x$ is Lebesgue-measurable. Let us then prove that for any $x \in [-n,n]$ and any $\delta >0$,
\begin{equation} \label{mesurabilite}
\{ z \in E \, : \, g(x, z) \geq \delta \} \in {\mathcal L}(E) \text{ and } \{ z \in E \, : \, g(x, z) \leq -\delta \} \in {\mathcal L}(E) 
\end{equation}
which will end the proof of the mesurability of $g$ and then of the lemma. For any $x \in [-n,n]$ and any $\delta >0$, we notice that
$$
\{ z \in E \, : \, g(x, z) \geq \delta \} =
  \begin{cases}
 \{ z\in E \, : \, H_z (x) >0\} \cap \{ z \in E \, : \, d(x, A_z) \geq \delta \}  \quad   
\text{if $|x| \geq \delta$}\\
  \ \emptyset \phantom{z \in E \, : \, H_z (x) >0\} \cap \{ z \in E \, : \, d(x, A_z) \geq \delta \}}  \ \text{ if $|x| < \delta$}.
  \end{cases}
$$
If $|x| \geq \delta$,
\begin{eqnarray*}
& & \{ z\in E \, : \, H_z(x) >0\} \cap \{ z \in E \, : \, d(x, A_z) \geq \delta \}  \\
& & \quad \quad = \{ z \in E \, : \, H_z (x) >0\} \cap \{ z\in E  \, : \, ]x - \delta, x + \delta [ \cap A_z = \emptyset \} \\
& & \quad \quad = \{ z \in E \, : \, H_z (y) >0 \, , \, \forall y \in  ]x - \delta, x + \delta [\} \, .
\end{eqnarray*}
Using the fact that the function $H_z (\cdot)$ is continuous, we get that if $|x| \geq \delta$,
\begin{eqnarray*}
& & \{ z \in E \, : \, H_z (x) >0\} \cap \{ z \in E \, : \, d(x, A_z) \geq \delta \}  \\
& & \quad \quad = \bigcup_{p \in {\mathbb N }^*} \big \{ z \in E \, : \, H_z(y) \geq \frac{1}{p} 
\, , \, \forall y \in  ]x - \delta, x + \delta [ \cap {\mathbb Q}\big \} \, ,
\end{eqnarray*}
which proves the first part of (\ref{mesurabilite}) since 
$\{ z \in E \, : \, H_z (a) \geq p^{-1} \}$ belongs to ${\mathcal L} (E)$ for any $a \in {\mathbb Q}$ and any $p \in {\mathbb N }^*$. The second part of (\ref{mesurabilite}) follows from the same arguments by changing the sign. This ends the proof of the ${\mathcal L} ( {\mathbb R} \times E) -{\mathcal B} ( {\mathbb R})$ measurability of  $ f$ defined by (\ref{defflemma}).

Next $P_{(U + N , Z)}$ and $P_{(V + N , Z)}$ are absolutely continuous wrt $\lambda \otimes P_Z$.  Consequently, starting from (\ref{lissage}) and using (\ref{formulemajw2}), the lemma follows.
$\diamond$
\begin{Proposition} \label{condversion} Let  $X_i = f(Y_i) - \bkE ( f(Y_i))$, where $f$ belongs to $\tMonm(Q, P_{Y_0})$. Assume that $M_{2,\alpha} (Q) < \infty$. Then the series
$\bkE (X_0^2) + 2 \sum_{k \geq 1} \bkE (X_0 X_k)$ is convergent to some nonnegative real  $\sigma^2$. If $\sigma^2 >0$, then there exists a positive constant 
$C$ depending on $\sigma^2$ such, that for any $n >0$,
\begin{equation}
\bkE  \big ( W_2^2 (P_{S_n | {\mathcal F}_0} , G_{n \sigma^2}) \big ) \leq C n^{1/2} 
M_{3,\alpha} ( Q , n^{1/2}) \, ,
\end{equation}
where $M_{3,\alpha} ( Q , n^{1/2})$ is defined in (\ref{3alphatronq}). 
\end{Proposition}

\noindent{\it Proof of Proposition \ref{condversion}.} 
Let  $\N$ be a sequence of independent random variables
with normal distribution ${\mathcal N}(0, \sigma^2)$. Suppose furthermore that the sequence $\N$ is independent of ${\mathcal F}_{\infty}$.
Let $N$ be a ${\mathcal N}(0,\sigma^2)$-distributed random variable, independent
of ${\mathcal F}_{\infty}\vee \sigma(N_i , i \in
{\mathbb Z}  )$. Set $T_n = N_1 + N_2 + \cdots + N_n$. Let $Z= ((Y_i, \delta_i) \, : \, i \leq 0)$ and $E = ({\mathbb R} \times [0,1] )^{{\mathbb Z}^-}$. Notice that 
$(E, {\mathcal L}(E), P_Z)$ is a purely non atomic Lebesgue space. From Lemma \ref{riolma}, we have to bound 
\begin{equation}\label{Delta}
\Delta (\varphi) =   {\mathbb E} ( \varphi(S_n + N, Z) - \varphi (T_n +  N,Z) )  \, ,
\end{equation}
for any function $\varphi$ in ${\Lambda_2 (E)}$. With this aim,  we apply the  Lindeberg method.

\begin{Notation}\label{not21}
Let $$\varphi_k (x,Z) = \int_{{\mathbb R}} \varphi (t, Z) \phi_{\sigma \sqrt{n-k+1}} (x-t) dt\, . $$ Let $S_0 = 0$, and, for $k>0$,
let $\Delta_k =  \varphi_k (S_{k-1} + X_k,Z) - \varphi_k ( S_{k-1} + N_k,Z)$.
\end{Notation}

Since the sequence $\N$ is independent of the sequence $\XZ$,
\begin{equation}\label{sumdelta}
{\mathbb E} ( \varphi (S_n +  N, Z) -\varphi (T_n +  N,Z) ) = \sum_{k=1}^n {\mathbb E} ( \Delta_k ) .
\end{equation}
We first show that for any real $u \in [0,1]$,
\begin{equation} \label{1resrio95}
|\bkE (\Delta_k)| \leq   C \big ( (n-k+1)^{-1/2} +  D_k(u) \big ) \, ,
\end{equation}
where
\begin{eqnarray} \label{2resrio95}
& & D_k(u)  = (n-k+1)^{1/2} \int_0^{\alpha (k)} Q(x)dx  + \sum_{ i > [k/2]} \int_0^{\alpha (i)}  Q^2(x)  dx \nonumber \\
& & + \int_0^{u} Q(x) R(x)  dx + (n-k+1)^{-1/2} \int_0^{1} Q (x) R(x) R(x\vee u)  dx
 \, .
\end{eqnarray} 
We now prove (\ref{1resrio95}). For the sake of brevity, write $ \varphi_k (x, Z) =  \varphi_k(x)$ and $ \varphi (x, Z) =  \varphi(x)$ (the derivatives are taken wrt $x$). By the Taylor formula at order 3,
$$
\big | \bkE \big (  \varphi_k (S_{k-1} + N_k) - \varphi_k ( S_{k-1}) - \frac{\sigma^2}{2} \varphi''_k ( S_{k-1}) \big ) \big | \leq
\frac{\| \varphi^{(3)}_k  \|_{\infty}}{6 } \bkE |N|^3 \, .
$$
Now Lemma 6.1 in Dedecker, Merlev\`ede and Rio (2009) gives that, almost surely,
\beq \label{cons1lmadmr}
\| \varphi^{(i)}_k\|_{\infty} \leq c_{i} \sigma^{2-i}(n-k+1)^{(2-i)/2} \mbox{ for any integer $i \geq 2$ } \, 
\eeq
where the $c_i$'s are universal constants.
Therefore
$$
\big | \bkE \big (  \varphi_k (S_{k-1} + N_k) - \varphi_k ( S_{k-1}) - \frac{\sigma^2}{2} \varphi''_k ( S_{k-1}) \big ) \big | \leq
C (n-k+1)^{-1/2} \, .
$$
Consequently to prove (\ref{1resrio95}), it remains to show that
\begin{equation} \label{3resrio95}
\big | \bkE \big (  \varphi_k (S_{k-1} + X_k) - \varphi_k ( S_{k-1}) - \frac{\sigma^2}{2} \varphi''_k ( S_{k-1}) \big ) \big |  \leq  C D_k(u)  \, ,
\end{equation}
where $D_k(u)$ is defined by (\ref{2resrio95}). To prove (\ref{3resrio95}), we follow the lines of the proof of Proposition 2(a) of Rio (1995-b) with $b_2 =\| \varphi^{(2)}_k\|_{\infty}$,  
$b_3 = \| \varphi^{(3)}_k\|_{\infty}$ and the modifications below. 
Since $f$ belongs to $\tMonm(Q, P_{Y_0})$, we can write
$$
X_i = \lim_{N \rightarrow \infty} {\mathbb L}^1 \sum_{\ell=1}^N a_{\ell,N}\big(  f_{\ell,N}(Y_i)
- \bkE(f_{\ell,N}(Y_i)) \big )\, ,
$$ 
with $f_{\ell,N}$
belonging to $\tMon( Q, P_{Y_0})$ and $\sum_{\ell=1}^N |a_{\ell,N}|
\leq 1$. For $u \in [0,1]$, let the function $g_u$ be defined by $g_u(x) = (x \wedge Q(u)) \vee (-Q(u))$. Since there exists a subsequence $m(N)$ tending to infinity such that $\sum_{\ell=1}^{m(N)} a_{\ell, m(N)} g_u \circ f_{\ell, m(N)}(Y_0)$ is convergent in $ {\mathbb L}^1$, for any $i \geq 0$, we define
\[
 \bar X_i=
 \lim_{N \rightarrow \infty} {\mathbb L}^1 \sum_{\ell=1}^{m(N)} a_{\ell, m(N)}\big( g_u \circ f_{\ell, m(N)}(Y_i)
- \bkE(g_u \circ f_{\ell,m(N)}(Y_i)) \big ) \quad \text{and} \quad  \tilde X_i=X_i - \bar X_i \, .
\]
Let also
$$
Q_u(x) := Q (x)\I_{x \leq u} \text{ and } \bar Q_u (x) := Q(x \vee u) \, .
$$
Since $Q_{|g_u \circ f_{\ell, m(N)} (Y_i)|} \leq \bar Q_u $, this means that $\bar X_i = r(Y_i) - \bkE (r(Y_i)) $ where $r$ belongs to $\tMonm(\bar Q_u, P_{Y_0})$.

By the Taylor integral formula,
\begin{eqnarray} \label{dt1}
 \varphi_k(  S_k )  - \varphi_k (S_{k-1} )   -  \varphi_k^\prime (S_{k-1})  X_k &  = &
X_k \int_0^1  (\varphi_k^\prime (S_{k-1} + v X_k  ) - \varphi_k^\prime (S_{k-1} )) dv   \nonumber \\
& =  & X_k \int_0^1  (\varphi_k^\prime (S_{k-1} + v X_k  ) -
       \varphi_k^\prime (S_{k-1} + v \bar X_k )) dv \nonumber \\ & + &
X_k\bar X_k \int_0^1 \!\!\! \int_0^1 v \varphi_k'' ( S_{k-1} + vv'\bar X_k) dv dv' .
\end{eqnarray}
The first term on right hand is bounded up by
$b_2 |X_k (X_k - \bar X_k)|/2$. Moreover
$$
\Big|
\int_0^1 \!\!\! \int_0^1 v \varphi_k'' ( S_{k-1} + vv'\bar X_k) dv dv'
-{1\over 2} \varphi_k'' (S_{k-1}) \Big| \leq {b_3\over 6} |\bar X_k|.
$$
Setting $h_u(x) = x- g_u(x)$, we get that for any $f$ belonging to $\tMon( Q, P_{Y_0})$,
\begin{eqnarray*}
& & \bkE \big | (f(Y_k)
- \bkE(f(Y_k)) ) ( h_u \circ f_{\ell}(Y_k)
- \bkE(h_u \circ f(Y_k) ))  \big | \\
& & \quad \quad \quad \quad  \leq  \bkE | f(Y_k)
h_u(f(Y_k))  | +3  \bkE  | f(Y_k) | \bkE |
h_u(f(Y_k)) | \, .
\end{eqnarray*}
Since $Q_{| f
(Y_k) |} \leq Q $ and $Q_{|h_u ( f (Y_k))|} \leq (Q - Q(u))_+ \leq Q_u $, we derive that
\begin{eqnarray*}
& & \bkE \big | (f(Y_k)
- \bkE(f(Y_k)) ) ( h_u \circ f_{\ell}(Y_k)
- \bkE(h_u \circ f(Y_k) ))  \big |\\
& & \quad \leq  \int_0^u Q^2(x) dx +3  \big( \int_0^1 Q(x)dx \big) \big( \int_0^u Q(x)dx \big) 
\leq 4 \int_0^u Q^2(x) dx \, ,
\end{eqnarray*}
by using Lemma 2.1(a) in Rio (2000). Now, by Fatou lemma,
\begin{eqnarray*}
& & \bkE | X_k ( X_k - \bar X_k) | \leq \liminf_{N \rightarrow \infty} \sum_{\ell  =1}^{m(N)} \sum_{j  =1}^{m(N)}
|a_{\ell, m(N)}| |a_{j,m(N)}|\\
& &  \times  \bkE \big | (f_{\ell, m(N)}(Y_k)
- \bkE(f_{\ell, m(N)}(Y_k)) ) ( h_u \circ f_{j, m(N)}(Y_k)
- \bkE(h_u \circ f_{j, m(N)}(Y_k) ))  \big | \, ,
\end{eqnarray*}
whence
\beq \label{dt2}
\bkE | X_k ( X_k - \bar X_k) |  \leq 4 \int_0^u Q^2(x) dx \, .
\eeq
Similarly using Lemma 2.1  in Rio (2000) and the fact that $Q_{| g_u\circ f
(Y_k) |} \leq \bar  Q_u$ for any $f$ belonging to $\tMon( Q, P_{Y_0})$, we derive that
\beq \label{dt2bis}
\bkE | X_k ( \bar X_k)^2 |  \leq  8 \int_0^{1} Q^2 (x) Q (x\vee u) dx  \, .
\eeq
It follows that
\begin{eqnarray} \label{dt3}
& & \big | \bkE (   \varphi_k( S_k )   - \varphi_k (S_{k-1} )    -  \varphi_k^\prime (S_{k-1}) X_k -
{1\over 2} \varphi_k'' (S_{k-1}) X_k \bar X_k ) \big |  \nonumber \\
& & \leq
2b_2 \int_0^u Q^2(x) dx +
{4b_3\over 3} \int_0^{1} Q^2 (x) Q (x\vee u) dx  .
\end{eqnarray}
Now we control the second order term.
Let \beq \label{defgamma} \Gamma_{k}(k,i) =  \varphi_k'' (S_{k-i}) - \varphi_k''  (S_{k-i-1}) \, ,\eeq
and
\beq \label{defr}
r = \alpha^{-1}(u) \, .
\eeq
Clearly
$$
\varphi_k'' (S_{k-1})  X_k \bar X_k = \sum_{i=1}^{(r\wedge k)-1}
\Gamma_{k}(k,i) X_k \bar X_k  + \varphi_k'' (S_{k-(r\wedge k)})X_k \bar X_k  \, ,
$$
Since
$|\Gamma_{k}(k,i) | \leq b_3 |X_{k-i}|$, by stationarity we get that for any $i \leq (r\wedge k)-1$,
$$
\big |\cov(\Gamma_{k}(k,i) , X_k \bar X_k  ) \big | \leq  b_3 \Vert X_0 \big (\bkE_0 (X_i \bar X_i )- \bkE (X_k \bar X_k) \big ) \Vert_1 \, .
$$
Applying Proposition \ref{covineg} with $m=1$, $q=2$, $k_1=0$, $k_2=k_3=i$, $f_{j_1}=f_{j_2}=f$ and
$f_{j_3} \in \tMonm( \bar Q_u, P_{Y_0})$, we derive that
\begin{eqnarray*}
& & \big |\cov (\Gamma_{k}(k,i) , X_k \bar X_k  ) \big | \leq 32 b_3 \int_0^{\alpha(i)} Q^2(x) Q (x \vee u) dx
  \, .
\end{eqnarray*}
Since $| \varphi_k'' (S_{k-(r\wedge k)}) | \leq b_2 $ a.s., we also get
by stationarity that
$$
\big |\cov (  \varphi_k'' (S_{k-(r\wedge k)}), X_k\bar X_k)) \big | \leq  b_2 \Vert \bkE_0 (X_{r\wedge k} \bar X_{r\wedge k})-  \bkE (X_{r\wedge k} \bar X_{r\wedge k})\Vert_1 \, .
$$
Applying Proposition \ref{covineg} with $m=0$, $q=2$, $k_1=k_2=r$, $f_{j_1}=f$ and
$f_{j_2} \in \tMonm(\bar Q_u, P_{Y_0})$, and noting that $\alpha (r) \leq u$, we also get that
\begin{eqnarray*}
\big |\cov (  \varphi_k'' (S_{k-(r\wedge k)}), X_k\bar X_k)) \big | &\leq  &16 b_2 \Big ( \int_0^{u} Q(x) Q (u) dx \I_{r \leq k} +
 \int_0^{{\alpha(k)}} Q(x) Q (x \vee u) dx \I_{k < r} \Big ) \, .
 \end{eqnarray*}
Hence
\begin{eqnarray*} \label{dt4}
 {1\over 2}  |\cov ( \varphi_k'' (S_{k-1}),
  X_k\bar X_k) | & \leq &  8 b_2 \int_0^{u} Q(x) Q ( u) dx\I_{r \leq k} +8 b_2 \int_0^{{\alpha(k)}} Q(x) Q (x \vee u) dx \I_{k < r}  \nonumber \\
& + &    16 b_3 \int_0^{1} Q^2(x) R (x \vee u) dx   \,  ,
\end{eqnarray*}
which together with (\ref{dt3}) and (\ref{dt2})  implies that
\begin{eqnarray} \label{dt5}
& & \big | \bkE (   \varphi_k( S_k )   - \varphi_k (S_{k-1} )    -  \varphi_k^\prime (S_{k-1}) X_k) -
{1\over 2} \bkE (\varphi_k'' (S_{k-1})) \bkE (X_k^2) \big | \leq
12b_2 \int_0^u Q^2(x) dx +
  \nonumber \\
& & 8 b_2 \int_0^{{\alpha(k)}} Q(x) Q (x \vee u) dx \I_{k < r}+ \frac{52}{3} b_3 \int_0^{1} Q^2(x) R (x \vee u) dx \,  .
\end{eqnarray}
To give now an estimate of the expectation of $\varphi'_k ( S_{k-1})X_k$, we write
$$
\varphi'_k (S_{k-1} ) = \varphi'_k (0) + \sum_{i=1}^{k-1} ( \varphi'_{k} (S_{k-i}) - \varphi'_{k} ( S_{k-i-1}) ) .
$$
Hence
\begin{eqnarray} \label{dt6}
\bkE (\varphi'_k ( S_{k-1})X_k) & = &  \sum_{i=1}^{k-1} {\rm Cov} \big (\varphi'_{k} (S_{k-i}) - \varphi'_{k} ( S_{k-i-1} ) , X_k \big ) +
\bkE ( \varphi'_k (0) X_k ) \, .
\end{eqnarray}
Now $\varphi'_k (0) $ is a ${\cal F}_0$-measurable random variable, and since $\varphi'(0)=0$ and $\varphi'$ is $1$-Lipschitz wrt $x$, 
$$
|\varphi'_k (0)| =  |\int (\varphi'(u) -\varphi'(0))  \phi_{\sigma \sqrt{ n-k+1 } } (-u) du | \leq  \sigma \sqrt{ n-k+1 } \, \text{ a.s.}
$$
Applying Proposition \ref{covineg} with $m=0$, $q=1$, $k_1=k$ and $f_{j_1}=f$, it follows that
\beq \label{dt7}
{\mathbb E} ( \varphi'_k (0) X_k ) \leq \sigma\sqrt{ n-k+1} \Vert \bkE_0 (X_k) \Vert_1 \leq 8 \sigma \sqrt{n-k+1} \int_0^{\alpha (k)/2} Q(x)dx \, .
\eeq
We give now  an estimate of $\sum_{i=1}^{k-1} {\rm Cov} \big (\varphi'_{k} (S_{k-i}) - \varphi'_{k} ( S_{k-i-1} ) , X_k \big )$. Using the 
stationarity and noting that $|\varphi_k^\prime (S_{k-i} ) -\varphi_k^\prime (S_{k-i-1})| \leq b_2|X_{k-i}|$,
we have
$$
|\cov (\varphi_k^\prime (S_{k-i} ) -\varphi_k^\prime (S_{k-i-1} ),  X_k)| \leq b_2 \Vert X_0 \bkE_0 ( X_i) \Vert_1 \, .
$$
Now, for any $ i \geq r$, $\alpha(i) \leq u$. So applying
Proposition \ref{covineg} with $m=1$, $q=1$, $k_1=0$, $k_2=i$,
$f_{j_1}=f_{j_2}=f$, we get, for any $k \geq i \geq r$, that
\begin{eqnarray} \label{dt8}
& & |\cov (\varphi_k^\prime (S_{k-i} ) -\varphi_k^\prime (S_{k-i-1} ),  X_k)| \leq 16 b_2 \int_0^{u} Q^2(x) \I_{x < \alpha (i) } dx  \, .
\end{eqnarray}
From now on, we assume that $i<r \wedge k$. Let us replace $X_k$ by $\bar X_k$. Since by stationarity,
$$
|\cov (\varphi_k^\prime (S_{k-i} ) -\varphi_k^\prime (S_{k-i-1} ),  X_k- \bar X_k)| \leq b_2 \Vert X_0 \bkE_0 ( X_i- \bar X_i) \Vert_1 \, ,
$$
we can apply Proposition \ref{covineg} with $m=1$, $q=1$, $k_1=0$, $k_2=i$, $f_{j_1}=f$ and
$f_{j_2} \in \tMonm( \bar Q_u, P_{Y_0})$. Consequently,
\begin{eqnarray} \label{dt9}
& & |\cov (\varphi_k^\prime (S_{k-i} ) -\varphi_k^\prime (S_{k-i-1} ),  X_k- \bar X_k)| \leq 16 b_2 \int_0^{u} Q^2(x) \I_{x < \alpha(i)}  dx  \, .
\end{eqnarray}
Now
$$
\varphi_k^\prime (S_{k-i} ) -\varphi_k^\prime (S_{k-i-1} ) - \varphi_k'' (S_{k-i-1})  X_{k-i}
=  R_{k,i},
$$
where $R_{k,i}$ is ${\cal F}_{k-i}$-measurable and
$|R_{k,i}| \leq b_3 X_{k-i}^2 /2$. Consequently, by stationarity,
$$
|\cov (R_{k,i}, \bar X_k) | \leq b_3 \Vert X^2_0 \bkE_0 ( \bar X_i) \Vert_1 /2 \, .
$$
Applying Proposition \ref{covineg} with $m=2$, $q=1$, $k_1=k_2=0$, $k_3=i$, $f_{j_1}=f_{j_2}=f$ and
$f_{j_3} \in \tMonm(\bar Q_u, P_{Y_0})$, we get that
\begin{eqnarray} \label{dt10}
|\cov (R_{k,i}, \bar X_k) |  & \leq & 32 b_3 \int_0^{ \alpha(i) } Q^2(x)Q(x \vee u)  dx 
 \, .
\end{eqnarray}
In order to estimate the term $\cov (\varphi''_k (S_{k-i-1})X_{k-i}, \bar X_k)$,
we introduce the decomposition below:
$$
\varphi''_k (S_{k-i-1}) = \sum_{l=1}^{(i-1) \wedge (k-i-1)} ( \varphi''_k
(S_{k-i-l}) - \varphi''_k (S_{k-i-l-1}) )  + \varphi''_k (S_{(k-2i)\vee 0}) .
$$
For any $l\in \{ 1, \cdots, (i-1) \wedge (k-i-1) \}$, by using the notation (\ref{defgamma}) and stationarity, we get that
$$
|\cov (\Gamma_{k}(k,l+i) X_{k-i} , \bar X_k )| \leq b_3 \Vert X_{- l} X_0 \bkE_0 (\bar X_i ) \Vert_1 \, .
$$
Applying Proposition \ref{covineg} with $m=2$, $q=1$, $k_1=-\ell$, $k_2=0$, $k_3=i$, $f_{j_1}=f_{j_2}=f$ and
$f_{j_3} \in \tMonm(\bar Q_u, P_{Y_0})$, we then derive that
\begin{eqnarray} \label{dt11}
|\cov (\Gamma_{k}( k, l+i)X_{k-i} , \bar X_k )|  & \leq & 64 b_3 \int_0^{ \alpha(i) } Q^2(x)Q(x \vee u)  dx  \, .
\end{eqnarray}
As a second step, we bound up
$|\cov ( \varphi_k'' (S_{(k-2i)\vee 0}) ,  X_{k-i}  \bar X_k )|$. Assume first that $i \leq [k/2]$. Clearly, using the notation (\ref{defgamma}),
$$
 \varphi_k'' (S_{k-2i}) = \sum_{l=i}^{(r-1) \wedge (k-i -1)} \Gamma_{k}(k, l+i) +  \varphi'' (S_{(k-i-r) \vee 0}).
$$
Now for any $l\in \{ i, \cdots, (r-1) \wedge (k-i-1) \}$, by  stationarity,
$$
|\cov (\Gamma_{k}(k, l+i)   , X_{k-i}\bar X_k )| \leq  b_3 \Vert X_{- l} \big (\bkE_{- l} (X_0\bar X_i )
- \bkE ( X_0 \bar X_i ) \big ) \Vert_1 \, .
$$
Hence applying Proposition \ref{covineg} with $m=1$, $q=2$, $k_1=-l$, $k_2=0$, $k_3=i$, $f_{j_1}=f_{j_2}=f$ and
$f_{j_3} \in \tMonm( \bar Q_u, P_{Y_0})$, we derive that
\begin{eqnarray} \label{dt12}
|\cov (\Gamma_{k}( k, l+i)  , X_{k-i}\bar X_k )| & \leq & 32 b_3 \int_0^{ \alpha(l) } Q^2(x)Q(x \vee u)  dx  \, .
\end{eqnarray}
If $i \leq k-r$, then stationarity implies that
$$
|\cov (\varphi_k'' (S_{k-i-r } ), X_{k-i}\bar X_k )|\leq  b_2 \Vert  \bkE_{0} (X_r\bar X_{i+r} )
- \bkE ( X_{r} \bar X_{i+r} ) \big ) \Vert_1 \, .
$$
Noting that $\alpha(r)\leq u < \alpha(i)$ and applying Proposition \ref{covineg} with $m=0$, $q=2$, $k_0=0$, $k_1=r$, $k_2=i+r$, $f_{j_1}=f$ and
$f_{j_2} \in \tMonm(\bar Q_u, P_{Y_0})$, we also get that
\beq \label{dt13}
|\cov (\varphi_k'' (S_{k-i-r } ), X_{k-i}\bar X_k )|
\leq 16b_2 \int_0^u \I_{x < \alpha(i)} Q(x) Q(u) dx  \, .
\eeq
Now if $i > k-r$, then we write that
$$
|\cov (\varphi_k'' (0 ), X_{k-i}\bar X_k )|\leq  b_2 \Vert  \bkE_{0} (X_{k-i}\bar X_k )
- \bkE (X_{k-i}\bar X_k )  \Vert_1 \, .
$$
Applying Proposition \ref{covineg} with $m=0$, $q=2$, $k_0=0$, $k_1=k-i$, $k_2=k$, $f_{j_1}=f$ and
$f_{j_2} \in \tMonm(\bar Q_u, P_{Y_0})$, and noting that for $i \leq [k/2]$, $ \alpha(k-i) \leq  \alpha([k/2])$, we obtain that
\beq \label{dt13pri}
|\cov (\varphi_k'' (0 ), X_{k-i}\bar X_k )|
\leq 16b_2 \int_0^{\alpha([k/2])} Q(x) Q(x \vee u) dx  \, .
\eeq
Assume now that $i \geq [k/2] + 1$. For any $i \leq k$, the stationarity entails that
$$
|\bkE ( \varphi_k''(0) X_{k-i} \bar X_k) |\leq b_2 \Vert X_0 \bkE_0
(\bar X_i) \Vert_1 \, .
$$
Hence applying Proposition \ref{covineg} with $m=1$, $q=1$, $k_0=0$, $k_1=i$, $f_{j_1}=f$ and
$f_{j_2} \in \tMonm(\bar Q_u, P_{Y_0})$, and noting that for $i \geq [k/2]+1$, $ \alpha(i) \leq  \alpha([k/2])$, we obtain that
\begin{eqnarray} \label{dt13second}
|\bkE ( \varphi_k''(0) X_{k-i} \bar X_k) |&\leq & 16 b_2 \int_0^{\alpha([k/2])}  Q(x) Q(x \vee u) dx   \, .
\end{eqnarray}
Adding the inequalities (\ref{dt7}), (\ref{dt8}),
(\ref{dt9}), (\ref{dt10}), (\ref{dt11}), (\ref{dt12}) (\ref{dt13}), (\ref{dt13pri}) and (\ref{dt13second}), summing on $i$ and $l$,  and using the fact that 
$$\sum_{i=1}^{k-1}\I_{x < \alpha (i)} \leq \alpha^{-1} (x) \, , \,  \sum_{i=1}^r  \I_{x < \alpha (i)}
\leq \alpha^{-1} (x \vee u)  \text{ and } \sum_{i=1}^r  i \I_{x < \alpha(i)} \leq
 (\alpha^{-1} (x \vee u))^2 \, , $$
we then get:
\begin{eqnarray} \label{dt14}
 & & \quad |\bkE (\varphi^\prime (  S_{k-1})  X_k )  -  \sum_{i=1}^{r-1}
\bkE (\varphi'' (S_{k-2i})) \bkE  (X_{k-i} \bar X_k ) \I_{i \leq [k/2]}| \leq  C (n-k+1)^{1/2} \int_0^{\alpha(k)} 
Q(x)dx +  \nonumber \\
 & & 
48 b_2 \int_0^{u} Q(x)R (x)  dx    + 24 k  b_2 \int_0^{\alpha([k/2])} Q(x) Q(x \vee u) dx   + 128 b_3 \int_0^{ 1 }   Q^2(x)R(x \vee u) dx 
\,  .
\end{eqnarray}
It remains to bound up
$$
A_k:=\sum_{i=1}^{r-1} \bkE (\varphi_k'' (S_{k-2i})) \bkE ( X_{k-i}
\bar X_k) \I_{i \leq [k/2]}- \sum_{i=1}^{\infty} \bkE (\varphi_k''
(S_{k-1}) )  \bkE (X_{k-i}X_k) \, .
$$
We first note that by stationarity,
$$
\sum_{i\geq r}| \bkE (\varphi_k'' (S_{k-1}) )  \bkE (X_{k-i}X_k) |
\leq  b_2 \sum_{i\geq r} | \bkE(f(Y_0) \bkE_0(X_i)) | \, .
$$
Applying Proposition \ref{covineg} and noting that $\alpha(i) \leq u$ for $i \geq r$, we get that
\begin{equation} \label{dt16}
\sum_{i\geq r}| \bkE (\varphi_k'' (S_{k-1}) )  \bkE (X_{k-i}X_k) | \leq  
8b_2 \sum_{i \geq r} \int_0^{\alpha (i)} Q^2(x)dx  \leq  8b_2 \int_0^u Q(x) R (x)dx \, .
\end{equation}
By stationarity we also have
$$
\sum_{i=1}^{r-1}| \bkE (\varphi_k'' (S_{k-1}) )  \bkE (X_{k-i}(X_k - \bar X_k)) |
\leq  b_2 \sum_{i=1}^{r-1} | \bkE(f(Y_0) \bkE_0(X_i-\bar X_i)) | \, .
$$
Next, noting that $u < \alpha (i)$ for all $i <r$ and
applying  Proposition \ref{covineg}, we get that
\begin{eqnarray} \label{dt17}
\sum_{i=1}^{r-1}| \bkE (\varphi_k'' (S_{k-1}) )  \bkE (X_{k-i}(X_k - \bar X_k)) | &\leq  &
8b_2 \int_0^u Q^2(x) \sum_{i =1}^{r-1}\I_{x <\alpha (i)}dx \nonumber \\
& \leq & 8b_2 \int_0^u Q^2(x)  \alpha^{-1} (x)dx \, .
\end{eqnarray}
In addition, another application of Proposition \ref{covineg} gives
\begin{eqnarray} \label{dt17bis}
\sum_{i=1+[k/2]}^{r-1}| \bkE (\varphi_k'' (S_{k-1}) )  \bkE (X_{k-i} \bar
X_k) |  &\leq  &
8b_2 \sum_{ i > [k/2]}  \int_0^{\alpha(i)} Q^2(x) dx  \, .
\end{eqnarray}
In order to bound up the last term, we still write
$$
\bkE (\varphi_k'' (S_{k-1}) - \varphi_k'' (S_{k-2i})) \bkE ( X_{k-i}
\bar X_k ) \I_{i \leq [k/2]} = \sum_{l=1}^{2i-1} \bkE (\Gamma_{k}(k,l))  \bkE (
f(Y_{0}) \bkE_0(\bar X_i) ) \I_{i \leq [k/2]}.
$$
Both this decomposition, Proposition \ref{covineg} and Lemma 2.1 in Rio (2000) then yield :
\begin{eqnarray} \label{dt18}
& & \sum_{i=1}^{r-1}| \bkE (\varphi_k'' (S_{k-1}) - \varphi_k'' (S_{k-2i})) \bkE ( X_{k-i} \bar X_k ) |\I_{i \leq [k/2]}
 \leq  8 b_3  \sum_{i =1}^{r-1} i \int_0^{\alpha(i)} Q^2(x) Q(x \vee u)dx \nonumber \\
&  & \quad \quad \quad \quad \leq 8 b_3  \int_0^1  Q(x)R(x) R(x \vee u) dx \, .
\end{eqnarray}
Hence (\ref{dt16}), (\ref{dt17}) and (\ref{dt18}) together entail that
\begin{equation} \label{dt19}
|A_k| \leq     16b_2 \int_0^u Q(x)  R (x)dx + 8b_2 \sum_{ i > [k/2]} \int_0^{\alpha(i)} Q^2(x)  dx + 8 b_3 \int_0^1 Q(x) R(x) R(x \vee u) dx \, .
\end{equation}
(\ref{dt19}), (\ref{dt14}), (\ref{dt5}) together with (\ref{cons1lmadmr}) then yield (\ref{1resrio95}).

Notice now that
$$
\sum_{k=1}^n\sqrt{n-k+1} \int_0^{\alpha(k)} Q(x)dx \leq n^{1/2}  \int_0^1 (\alpha^{-1} (x)\wedge n ) Q(x) dx \, ,
$$
and that
\begin{eqnarray*}
& & \sum_{k=1}^n\sum_{ i > [k/2]}  \int_0^{\alpha(i)}  Q^2  (x) dx  \leq  2 \sum_{ i \geq 1 } (i \wedge n)
\int_0^{\alpha(i)} Q^2(x)dx \\
&  & \quad \leq 2
\int_0^1 Q(x) R(x) (\alpha^{-1}(x)\wedge n )dx 
\leq 2n^{1/2} \int_0^1 Q(x)R(x) (R(x)\wedge n^{1/2} )dx \, .
\end{eqnarray*}
Moreover 
\begin{eqnarray*}
n^{1/2}  \int_0^1 (\alpha^{-1} (x)\wedge n ) Q(x) dx  \leq  n^{1/2}\int_0^1 Q(x)R(x) (R(x)\wedge n^{1/2} )dx \, .
\end{eqnarray*}
\par
Hence to prove Proposition \ref{condversion}, it remains to select $u=u_k$ in such a way that
\begin{equation} \label{b2}
\sum_{k=1}^n \int_0^{u_k} Q (x) R ( x) dx  + \sum_{k=1}^n 
\frac{1}{\sqrt{k}} \int_0^1 Q(x)R(x) R(x\vee u_k) dx
 \leq C n^{1/2} M_{3,\alpha} (Q ,n^{1/2})\, .
\end{equation}
Let $R^{-1}(y) = \inf\{ v \in [0,1] \, : \, R(v) \leq y \}$  be the right continuous inverse of $R$. 
Since $R$ is right continuous, $x < R^{-1}(y)$  if and only if $R(x) > y$. We now choose 
$u_k = R^{-1}( k^{1/2} )$, so that 
\begin{equation}\label{defuk}
R (u_k) \leq k^{1/2} \ \text{and}\ R(x) > k^{1/2} \ \text{for any}\ x<u_k.
\end{equation}
With this choice of $u_k$, on one hand,
\begin{eqnarray} \label{b3}
  \sum_{k=1}^n \int_0^{u_k} Q(x) R(x) dx  & =  &
 \int_0^{1} Q(x) R(x) \sum_{k=1}^n \I_{R(x) > \sqrt{k}} dx  \leq \int_0^{1} Q(x)R(x)  (R^2(x) \wedge n) dx \nonumber \\
& \leq & n^{1/2}\int_0^{1}  Q (x) R(x)  (R(x) \wedge n^{1/2}) dx \, .
\end{eqnarray}
On the other hand
\begin{equation} \label{b40}
 \sum_{k=1}^n  \frac{1}{\sqrt{k}} \int_0^{1} Q(x) R(x)R (x\vee u_k) dx 
 \leq  
\sum_{k=1}^n \frac{1}{\sqrt{k}} \int_{u_k}^1 Q(x) R^2 (x) dx \\
 +\sum_{k=1}^n \int_0^{u_k} Q(x) R(x) dx  
  \end{equation}
using (\ref{defuk}). Next 
\begin{equation} \label{b42}
\sum_{k=1}^n \frac{1}{\sqrt{k}} \int_{u_k}^1 Q(x)R^2 (x) dx  \leq 
\sum_{k=1}^n \frac{1}{\sqrt{k}} \int_{u_n}^1 Q(x)R^2 (x) dx 
 \leq  2 n^{1/2} M_{3,\alpha} (Q ,n^{1/2})\,  .
\end{equation}
Combining (\ref{b40}) with (\ref{b42}) and (\ref{b3}), we then get  (\ref{b2}) ending the proof of the proposition.
$\diamond$

\begin{Proposition} \label{inegamax} For $f$ in $\tMonm( Q, P_{Y_0})$,
let  $X_i = f(Y_i) - \bkE ( f(Y_i))$. Set $S_n^* = \max_{1 \leq k\leq n}| S_k|$. Assume that $M_{2,\alpha} (Q) < \infty$. Then the series
$\bkE (X_0^2) + 2 \sum_{k \geq 1} \bkE (X_0 X_k)$ is convergent to some nonnegative real  $\sigma^2$ and for any positive real $\lambda$,
\begin{eqnarray*}
\BBP (S_n^* \geq 5 \lambda) & \leq & c_1\exp \Bigl( - \frac{\lambda^2}{c_2n
\sigma^2} \Bigr) + c_3 n \lambda^{-3}  M_{3,\alpha} (Q , \lambda) + c_4 n \sigma^3 \lambda^{-3} \, ,
\end{eqnarray*}
where $M_{3,\alpha} ( Q , n^{1/2})$ is defined in (\ref{3alphatronq}) and $c_1$, $c_2$, $c_3$ and $c_4$ are positive constants not depending on $\sigma^2$, so that the first term vanishes if $\sigma^2 =0$.
\end{Proposition}

\noindent{\bf Proof of Proposition \ref{inegamax}.}
Assume first that $X_i = \sum_{\ell=1}^L a_{\ell} f_{\ell}(Y_i)
- \sum_{\ell=1}^L a_{\ell}\bkE(f_{\ell}(Y_i))$, with $f_\ell$
belonging to $\tMon(Q, P_{Y_0})$ and $\sum_{\ell=1}^L |a_\ell|
\leq 1$. Let $M>0$ and $g_M(x) = (x \wedge M) \vee (-M)$. For
any $i \geq 0$, we first define
\[
  X_i'=
 \sum_{\ell=1}^L a_{\ell}\big( g_M \circ f_{\ell}(Y_i)
- \bkE(g_M \circ f_{\ell}(Y_i)) \big ) \quad \text{and} \quad  X_i''=X_i -  X_i' \, .
\]
Let $q$ be a positive integer such that $q \leq n$. Let us first show that
  \begin{equation} \label{dec1FN}
  \max_{1 \leq k \leq n } |S_k | \leq
  \max_{ 1 \leq k \leq n}  | \bkE (S_n |\F_k)|
  + 2 q M +   \max_{1 \leq k \leq n }  \bkE_k \big (\sum_{i=1}^n |X_i''| \big ) + \max_{1 \leq k \leq n } \bkE_k \big (\sum_{i=1}^n |\bkE_{i-q} (X_i')| \big ) \, .
  \end{equation}
Notice that
$$
S_k = \bkE (S_n |\F_k) - \sum_{i=k+1}^n \bkE (X''_i |\F_k) - \sum_{i=k+1}^n \bkE (X'_i |\F_k) \, .
$$
Now
$$
\sum_{i=k+1}^n \bkE (X'_i |\F_k)  =  \sum_{i=k+1}^n \bkE (X'_i - \bkE_{i-q} (X_i') |\F_k)  - \sum_{i=k+1}^n \bkE ( \bkE_{i-q} (X_i') |\F_k) \, .
$$
The inequality (\ref{dec1FN}) follows by noticing that
$$
 \sum_{i=k+1}^n \bkE (X'_i - \bkE_{i-q} (X_i') |\F_k)  = \sum_{i=k+1}^{q+k} (\bkE_k (X'_i) - \bkE_{i-q} (X_i') )  \leq 2qM \, .
$$
Notice now that $(\bkE (S_n |\F_k))_{k \geq 1}$, $\Big (\bkE_k \big (\sum_{i=1}^n |X_i''| \big ) \Big )_{k \geq 1}$ and
$\Big (\bkE_k \big (\sum_{i=1}^n |\bkE_{i-q} (X_i')| \big ) \Big )_{k \geq 1}$ are martingales with respect to the filtration
$(\F_k)_{k \geq 1}$. Consequently from (\ref{dec1FN}) and the Doob maximal inequality, we infer that for any nondecreasing, non negative, convex and even function
$\varphi$ and if $qM \leq \lambda$,
  \begin{eqnarray} \label{decproba}
 \BBP (S_n^* \geq 5 \lambda) &\leq & \frac{\bkE(\varphi(S_n))}{\varphi(\lambda)} + \lambda^{-1}\sum_{i=1}^n \bkE |X_i''| +
\lambda^{-1}\sum_{i=1}^n  \| \bkE_{i-q} (X_i') \|_1 \, .
  \end{eqnarray}
Choose $u = R^{-1}(\lambda)$, $q = \alpha^{-1} (u) \wedge n $ and $M = Q (u)$. Since $R$ is right
continuous, we have $R(u)\leq \lambda$, hence $qM \leq R(u)  \leq \lambda$. Note also that 
  \begin{equation} \label{dec13FN}
  \sum_{k=1}^n  \bkE ( | X_k''|)
  \leq  2 n \int_0^u  Q (x) dx
  \leq  2n \int_0^1 Q(x) \I_{R(x)> \lambda} dx \, .
  \end{equation}
In addition using Proposition \ref{covineg}, we get that
\begin{equation}
\label{majnorm1prime}
\| \bkE ( X'_i | {\mathcal F}_{i-q} )\|_1 \leq 8 \int_0^{\alpha (q)/2} Q(x) dx  \, .
\end{equation}
Since $\alpha(q) /2 \leq u$,
$$
\sum_{i=1}^n  \| \bkE_{i-q} (X_i') \|_1 \leq 8 n  \int_0^1 Q(x) \I_{R(x) > \lambda} dx \, .
$$
It follows that
\begin{eqnarray} \label{decsecprime}
\lambda^{-1} \big ( \sum_{i=1}^n \bkE |X_i''| +
\sum_{i=1}^n  \| \bkE_{i-q} (X_i') \|_1 \big )
 & \leq & 10n \lambda^{-1} \int_0^1 Q(x) \I_{R(x) > \lambda} dx \nonumber \\ & \leq & 10n \lambda^{-2} \int_0^1 Q(x) R(x)  \I_{R(x) > \lambda} dx\, .
  \end{eqnarray}
To control now the first term in the inequality (\ref{decproba}), we choose the even convex function $\varphi$ such that
$$
\varphi(t) =  \left\{
\begin{array}{ll}0 & \text{if $0 \leq t \leq \lambda/2$} \\
\frac{1}{6} ( t - \frac{\lambda}{2} )^3 & \text{if $\lambda/2 \leq t \leq \lambda $} \\
\frac{\lambda^3}{48} + \frac{\lambda}{4} ( t - \lambda)^2+ \frac{\lambda^2}{8} ( t - \lambda) & \text{if $t \geq \lambda $} \, .
\end{array}
\right.
$$
Clearly $ \|\varphi^{(2)}\|_{\infty} \leq \lambda/2$ and  $\|\varphi^{(3)}\|_{\infty} \leq 1$. Let  $\N$
be a sequence of independent random variables
with normal distribution ${\mathcal N}(0, \sigma^2)$. Suppose furthermore that the sequence $\N$ is independent of $\X$.
 Set $T_n = N_1 + N_2 + \cdots + N_n$ and
 $\varphi_k (x) = {\mathbb E} ( \varphi (x + T_n - T_k)  )$.  With this notation
$$
{\mathbb E} ( \varphi (S_n) - \varphi (T_n) ) = \sum_{k=1}^n {\mathbb E} (  \varphi_k (S_{k-1} + X_k) - \varphi_k ( S_{k-1} + Y_k) ).
$$
To bound up ${\mathbb E} (  \varphi_k (S_{k-1} + X_k) - \varphi_k ( S_{k-1} + Y_k) ) $, we proceed as in the proof of Proposition \ref{condversion} with the following modifications. Firstly, $b_2 = \| \varphi^{(2)}_k\|_{\infty} \leq \lambda/2 $
and $b_3 = \| \varphi^{(3)}_k\|_{\infty} \leq 1 $. The following convention is also used: $S_0 = 0$ and for any positive integer $j$, $S_{-j}= - \sum_{i=1}^{j}X_{1-i}$. Notice that here the $\varphi_k$ are deterministic. Consequently $\bkE ( \varphi'_k (0) X_k) = 0$ and $\varphi_k'' (S_{\ell})$ is always 
${\mathcal F}_{\ell}$-measurable for any $\ell \in {\mathbb Z}$. We then infer that the following bound is valid:  for any $k=1, \dots, n$,
$$
{\mathbb E} (  \varphi_k (S_{k-1} + X_k) - \varphi_k ( S_{k-1} + Y_k) )  \leq   
\sigma^3 + C \lambda \int_0^u  Q(x) R(x) dx + C\int_0^1  Q (x) R(x) R(x\vee u) dx  \, ,
$$
where $C$ is a positive constant not depending on $\sigma^2$.
Choosing $u= R^{-1}(\lambda)$, we get that
\begin{eqnarray*}
\int_0^{u} Q(x) R(x) dx  = \int_0^{1} Q(x) R(x)  \I_{R(x) >  \lambda} dx \, ,
\end{eqnarray*}
and
\begin{eqnarray*}
\int_0^{1} Q(x)R(x) R(x\vee u) dx  \leq  \int_0^{1}  Q (x)  R(x) (R(x) \wedge \lambda)  dx \, .
\end{eqnarray*}
It follows that
\begin{equation} \label{1resphi}
{\mathbb E} ( \varphi (S_n) - \varphi (T_n) ) \leq  n \sigma^3 + 2 C n 
M_{3,\alpha} (Q, \lambda) \, .
\end{equation}
It remains to compute ${\mathbb E} (  \varphi (T_n) ) $. We have that
$6{\mathbb E} (  \varphi (T_n) ) \leq  {\mathbb E} \big (  T_n - \lambda/2 \big)_+^3$.
Hence, using the fact that $t^2 = \lambda^2/4 + (t-\lambda/2)^2 + \lambda (t-\lambda/2)$, we
obtain:
$$\bkE(\varphi(T_n)) \leq \frac{e^{-\lambda^2/(8n\sigma^2)}}{6}
\int_{0}^{\infty} e^{-\lambda x/(2n \sigma^2)} \frac{x^3}{\sigma \sqrt{2 n
\pi}} dx \, .
$$
Using the change of variables $y= \lambda x/(2n\sigma^2)$,  we derive that
\begin{equation} \label{majgauss} \bkE(\varphi(T_n)) \leq  \frac{\lambda^3}{\sqrt{2\pi} }
\Bigl( \frac{(2n\sigma^2)}{\lambda^2} \Bigr)^{7/2} e^{-\lambda^2/(8n\sigma^2)} \, .
\end{equation}
Starting from (\ref{decproba}) and collecting the bounds (\ref{decsecprime}), (\ref{1resphi}) and (\ref{majgauss}), the proposition is proved for any variable $X_i=f(Y_i) - \bkE(f(Y_i)) $ with
$f=\sum_{\ell=1}^L a_{\ell} f_{\ell}$ and $f_{\ell} \in \tMon (Q, P_{Y_0})$,
$\sum|a_\ell|\leq 1$. Since these functions are dense in $\tMonm (Q, P_{Y_0})$ by
definition, the result follows by applying Fatou's lemma. 

\begin{Proposition} \label{covineg} Let $m$ and $q$ be two nonnegative integers. 
For any $(m+q)$-tuple of integers $(j_{\ell})_{1\leq  \ell \leq m+q}$, let $X^{(j_{\ell})}_i = f_{j_{\ell}}(Y_i) - \bkE ( f_{j_{\ell}}(Y_i))$, where  $f_{j_{\ell}}$ belongs to $\tMonm( Q_{j_{\ell}}, P_{Y_0})$ for $1 \leq \ell \leq m+q$. Suppose that $Q_{j_{\ell}}^q$ 
is integrable for $\ell \geq m+1$.  Define the coefficients
 $\alpha_{k, {\bf Y}}(n) $ as in
(\ref{defalpha}).  Then for any integers $(j_{\ell})_{1 \leq \ell \leq m+q}$ and any integers $(k_{\ell})_{1 \leq \ell \leq m+q}$ such that
$k_1\leq k_2 \leq \cdots \leq k_{m+q}$ and $k_{m+1}-k_m = \ell$,  \begin{eqnarray*} \label{cov1}
\left \| \prod_{i=1}^{m}X^{(j_{i})}_{k_i}\Big (  \bkE_{k_m} \big ( \prod_{i=m+1}^{m+q}X^{(j_i)}_{k_i} \big ) - \bkE \big (  \prod_{i=m+1}^{m+q}X^{(j_i)}_{k_i}  \big ) \Big )
 \right \|_1 \leq 2^{m+q+2}  \int_0^{2^{q-2} \alpha_{q, {\bf Y}}(\ell) } \prod_{i=1}^{m+q} Q_{j_i}(x) dx \,  ,
 \end{eqnarray*}
 and
 \begin{eqnarray*} \label{cov1}
\left \| \prod_{i=1}^{m}f_{j_{i}}(Y_{k_i})\Big (  \bkE_{k_m} \big ( \prod_{i=m+1}^{m+q}X^{(j_i)}_{k_i} \big ) - \bkE \big (  \prod_{i=m+1}^{m+q}X^{(j_i)}_{k_i}  \big ) \Big )
 \right \|_1 \leq 2^{q+2}  \int_0^{2^{q-2} \alpha_{q, {\bf Y}}(\ell) }\prod_{i=1}^{m+q} Q_{j_i}(x) dx \,  ,
 \end{eqnarray*}
with the convention that $\prod_{i=1}^0 =\prod_{i=m+1}^m = 1$.
\end{Proposition}
{\bf Proof of proposition \ref{covineg}.} Assume first that  $f_{j_{\ell}}= \sum_{r=1}^N a_{r} g_{j_{\ell},r}$ where $ \sum_{r=1}^N |a_{r}| \leq 1$ and $g_{j_{\ell},r}$ belongs to $\tMon( Q_{j_{\ell}}, P_{Y_0})$ for $1 \leq \ell \leq m+q$. To soothe the notation, let also
\beq \label{notasimpX}
X^{(j_{\ell})}_{i, r} = g_{j_{\ell},r}(Y_i)
- \bkE(g_{j_{\ell}, r }(Y_i))  \, .
\eeq
We then have that
\begin{eqnarray*}
& & \left \| \prod_{i=1}^{m}X^{(j_i)}_{k_i}\Big (  \bkE_{k_m} \big ( \prod_{i=m+1}^{m+q}X^{(j_i)}_{k_i} \big ) - \bkE \big (  \prod_{i=m+1}^{m+q}X^{(j_i)}_{k_i}  \big ) \Big )
 \right \|_1 \\
 & & \leq \prod_{p=1}^{m+q} \big ( \sum_{r_p  =1}^N
|a_{r_p}| \big )
\left \| \prod_{i=1}^{m}X^{(j_i)}_{k_i, r_i}\Big (  \bkE_{k_m} \big ( \prod_{i=m+1}^{m+q}X^{(j_i)}_{k_i, r_i} \big ) - \bkE \big ( \prod_{i=m+1}^{m+q}X^{(j_i)}_{k_i, r_i} \big ) \Big )
 \right \|_1  \, .
\end{eqnarray*}
Now setting
$$A: =\Big | \prod_{i=1}^{m}X^{(j_i)}_{k_i, r_i} \Big | \sign \left \{ \bkE_{k_m} \big ( \prod_{i=m+1}^{m+q}X^{(j_i)}_{k_i, r_i} \big ) - \bkE \big ( \prod_{i=m+1}^{m+q}X^{(j_i)}_{k_i, r_i} \big )\right \} , $$ we get that
\begin{eqnarray*}
& & \left \|\prod_{i=1}^{m}X^{(j_i)}_{k_i, r_i}\Big (   \bkE_{k_m} \big ( \prod_{i=m+1}^{m+q}X^{(j_i)}_{k_i, r_i} \big ) - \bkE \big ( \prod_{i=m+1}^{m+q}X^{(j_i)}_{k_i, r_i} \big )\Big )
 \right \|_1 \\
& &   =  \bkE \left ( A  \Big (\bkE_{k_m} \big ( \prod_{i=m+1}^{m+q}X^{(j_i)}_{k_i, r_i} \big ) - \bkE \big ( \prod_{i=m+1}^{m+q}X^{(j_i)}_{k_i, r_i} \big ) \Big )\right ) 
  =  \bkE \left ( (A - \bkE(A))   \prod_{i=m+1}^{m+q}X^{(j_i)}_{k_i, r_i} \right )  \, .
  \end{eqnarray*}
From Proposition A.1 and Lemma A.1 in Dedecker and Rio (2008), we have that
\begin{eqnarray*}
\bkE \left ( (A - \bkE(A))   \prod_{i=m+1}^{m+q}X^{(j_i)}_{k_i, r_i}   \right )    \leq 2^{q+2}  \int_0^{\bar \alpha /2} Q_{|A|}(x) \prod_{i=m+1}^{m+q} Q_{j_i}(x)dx\, ,
\end{eqnarray*}
where %
\[
  \bar \alpha = \sup_{(t_1, \ldots , t_{q+1}) \in {\mathbf R}^{q+1}}\Big | \bkE \big ( (\I_{ A
  \leq t_1} - \p (A \leq t_1)) \prod_{i=m+1}^{m+q}(\I_{ g_{j_i,r_i} (Y_{k_i}) \leq t_{i- m+1} } - 
\p ( g_{j_i,r_i} (Y_{k_i}) \leq t_{i - m+1}) ) \big )\big | \, .\]
By monotonocity of the functions $g_{j_i,r_i}$, we then get that
\begin{eqnarray*}
\bar \alpha &\leq & 2^{q}\sup_{(t_1, \ldots ,t_{q+1})  \in {\mathbf R}^{q+1}} \big | \bkE \big ( (\I_{ A
  \leq t_1} - \p (A \leq t_1)) \prod_{i=m+1}^{m+q}(\I_{Y_{k_i} \leq t_{i- m+1} } - \p (Y_{k_i} \leq t_{i - m+1})) \big ) \big | \\ & \leq &  2^{q-1} \alpha_{q, {\bf Y}}(\ell) \, .
\end{eqnarray*}
Consequently,
\begin{eqnarray*}
& & \bkE \left ( (A - \bkE(A))   \prod_{i=m+1}^{m+q}X^{(j_i)}_{k_i, r_i}   \right )     \leq   2^{q+2}  \int_0^{2^{q-2} \alpha_{q, {\bf Y}}(\ell) } Q_{|A|}(x) \prod_{i=m+1}^{m+q} Q_{j_i}(x)dx \\
& & \leq  2^{q+2}  \int_0^{2^{q-2} \alpha_{q, {\bf Y}}(\ell) } \prod_{i=1}^m \big ( Q_{j_i}(x) + \int_0^1 Q_{j_i}(x) dx \big ) \prod_{i=m+1}^{m+q} Q_{j_i}(x)dx \, .
\end{eqnarray*}
Hence taking into account that $\prod_{i=1}^{m+q} \big ( \sum_{r_i  =1}^N
|a_{r_i}| \big )   \leq 1$ and using Lemma 2.1  in Rio (2000), the inequality is proved for functions $f_{j_{\ell}}= \sum_{r=1}^N a_{r} g_{j_{\ell},r}$ where $ \sum_{r=1}^N |a_{r}| \leq 1$ and $g_{j_{\ell},r}$ belongs to $\tMon( Q_{j_{\ell}}, P_{Y_0})$ for $1 \leq \ell \leq m+q$.

It remains to prove that the inequality remains valid for $f_{j_{\ell}}$ belonging to 
$\tMonm( Q_{j_{\ell}}, P_{Y_0})$ for $1 \leq \ell \leq m+q$. By definition, 
$$
X^{(j_{\ell})}_{i} = \lim_{N \rightarrow \infty} {\mathbb L}^1 \sum_{r=1}^N a_{r,N}X^{(j_{\ell})}_{i, r, N}   \, ,
$$
where $\sum_{r=1}^N |a_{r,N} |\leq 1$ and $
X^{(j_{\ell})}_{i, r, N} = g_{j_{\ell},r,N}(Y_i)
- \bkE(g_{j_{\ell}, r, N }(Y_i))  $ with the $g_{j_{\ell}, r, N }$ belonging to 
$\tMon( Q_{j_{\ell}}, P_{Y_0})$ for $1 \leq \ell \leq m+q$. 
Hence, by Fatou lemma the proposition will hold  if we can prove that the 
following inequality holds almost surely
\begin{eqnarray} \label{argfatou}
& & \bkE_{k_m} \big ( \prod_{i=m+1}^{m+q}X^{(j_i)}_{k_i} \big ) - \bkE \big (  \prod_{i=m+1}^{m+q}X^{(j_i)}_{k_i} \big )
\nonumber \\
& & \quad \quad =\lim_{N \rightarrow \infty} \prod_{i=m+1}^{m+q} \big (\sum_{r_i=1}^N a_{r_i,N} \big )
\Big (\bkE_{k_m} \big ( \prod_{i=m+1}^{m+q}X^{(j_i)}_{k_i, r_i, N} \big ) - \bkE \big (  \prod_{i=m+1}^{m+q}X^{(j_i)}_{k_i, r_i, N} \big )  \Big ) \, .
\end{eqnarray}
With this aim, notice that for any $m+1 \leq \ell \leq m+q$,
$$
X^{(j_{\ell})}_{i} =  \sum_{r=1}^N a_{r,N}X^{(j_{\ell})}_{i, r, N}   + \epsilon^{(j_{\ell})}_{i,N} \, ,
$$
with $\lim_{N \rightarrow \infty} \Vert \epsilon^{(j_{\ell})}_{i,N} \Vert_1 =0$. In addition, since for $m+1 \leq \ell \leq m+q$, $Q_{j_\ell}^q$ is integrable and $g_{j_{\ell}, r, N }$ belongs to $\tMon( Q_{j_{\ell}}, P_{Y_0})$, 
it follows that $\Vert X^{(j_{\ell})}_{i, r, N} \Vert_q \leq 2 \Vert Q_{j_\ell}\Vert_q$ and next 
$\Vert X^{(j_{\ell})}_{i}\Vert_q \leq 2 \Vert Q_{j_\ell}\Vert_q$ by an application of Fatou lemma. Consequently the $\epsilon^{(j_{\ell})}_{i,N}$'s are in ${\mathbb L}^q$ and satisfy 
$\Vert \epsilon^{(j_{\ell})}_{i,N} \Vert_q \leq 4 \Vert Q_{j_\ell}\Vert_q$. Now 
\begin{eqnarray*}
\Vert  \epsilon^{(j_{m+1})}_{k_{m+1},N} \prod_{i=m+2}^{m+q}X^{(j_i)}_{k_i} \Vert_1  & \leq & 2^{q-1}
\int_0^1 Q_{|\epsilon^{(j_{m+1})}_{k_{m+1},N}|} (x)   \prod_{i=m+2}^{m+q}Q_{j_i}(x) dx  \\
& \leq & 2^{q-1}
\int_0^1 Q_{|\epsilon^{(j_{m+1})}_{k_{m+1},N}|} (x)  Q^{q-1}_*(x) dx \,
 ,
  \end{eqnarray*}
where $Q_*= \max_{m+2 \leq i \leq m+q}Q_{j_i}$. Now for any positive $M$, $Q^{q-1}_* \leq M^{q-1} + Q^{q-1}_* \I_{Q_* >M}$. Hence,
\begin{eqnarray*}
2^{1-q} \Vert  \epsilon^{(j_{m+1})}_{k_{m+1},N} \prod_{i=m+2}^{m+q}X^{(j_i)}_{k_i} \Vert_1  & \leq & M^{q-1} \Vert \epsilon^{(j_{m+1})}_{k_{m+1},N}\Vert_1 +  \Vert \epsilon^{(j_{m+1})}_{k_{m+1},N} \Vert_q \Vert  Q_* \I_{Q_* >M} \Vert_q^{q-1}  \\
 & \leq & M^{q-1} \Vert \epsilon^{(j_{m+1})}_{k_{m+1},N}\Vert_1 + 4 \Vert Q_{j_{m+1}}\Vert_q  \Vert  Q_* \I_{Q_* >M} \Vert_q^{q-1} \, ,
 \end{eqnarray*}
 which tends to zero by letting first $N$ tends to infinity and after $M$. Similarly, we can show that for any $\ell \in \{ 1, \dots, q-1 \}$, 
 $$
 \lim_{N \rightarrow \infty}\Big \|  \Big ( \prod_{i=1}^{\ell} X^{(j_{m+i})}_{k_{m+i}, r, N} \Big ) \epsilon^{(j_{m+\ell +1})}_{k_{m+\ell + 1},N} \prod_{i=m+\ell +2}^{m+q}X^{(j_i)}_{k_i} \Big \|_1 =0 \, .
 $$
This ends the proof of (\ref{argfatou}) and then of the proposition. 
 

\end{document}